\documentclass[twoside]{article}
\usepackage[a4paper]{geometry}
\usepackage[latin1]{inputenc} % ou \usepackage[utf8]{inputenc}
\usepackage[T1]{fontenc} % ou \usepackage[OT1]{fontenc}
\usepackage{hyperref}
\usepackage{graphicx, amsmath, theorem, amssymb}

\usepackage{algorithm}
\usepackage{algorithmic}

\usepackage{color}
\usepackage{ wasysym }

\newtheorem{theorem}{Theorem}
\newtheorem{pr}{Proposition}

\newtheorem{Rm}{Remark}

\begin{document}

\newcommand{\blue}{\color{black}}

\newcommand{\proof}{{\bf Proof:~}}
\newcommand{\uio}{$\Sigma^{(0)}$}
\newcommand{\uioo}{$\Sigma^{(1)}$}
\newcommand{\uiok}{$\Sigma^{(k)}$}
\newcommand{\eorc}{$EORC$}
\newcommand{\uiokm}{$\Sigma^{(\overline{k})}$}
\newcommand{\uioK}{$\Sigma^{(k+1)}$}

\title{Rank Conditions for Observability and Controllability for Time-varying Nonlinear Systems}% sensor fusion}

\author{Agostino Martinelli\\INRIA, Grenoble, France\\agostino.martinelli@inria.fr}

%\address{INRIA, Grenoble, France}  % Please supply                                              

\maketitle

\begin{abstract}
This paper provides the extension of the observability rank condition and the extension of the controllability rank condition to time-varying nonlinear systems.
Previous conditions to check the state observability and controllability, only account for nonlinear systems that do not explicitly depend on time, or, for time-varying systems, they only account for the linear case.
In this paper, the general analytic conditions are provided. The paper shows that both these two new conditions (the extended observability rank condition and the extended controllability rank condition) reduce to the well known rank conditions for observability and controllability in the two simpler cases of time-varying linear systems and  time-invariant nonlinear systems.
{\blue
The proposed new conditions work automatically and can deal with any system, independently of its complexity (state dimension, type of nonlinearity, etc).
Simple examples illustrate both these conditions. 
In addition, the two new conditions are used to study the observability and the controllability properties of a lunar module. For this system, the dynamics exhibit an explicit time-dependence due to the variation of the weight and the variation of the moment of inertia. These variations are a consequence of the fuel consumption. To study the observability and the controllability properties of this system, the extended observability rank condition and the extended controllability rank condition introduced by this paper are required.
The paper shows that, even under the constraint that the main rocket engine delivers constant power, the state is weakly locally controllable. Additionally, it is weakly locally observable up to the yaw angle.

}
\end{abstract}

\newpage

\tableofcontents

\newpage

%%%%%%%%%%%%%%%%%%%%%%%%%%%%%%%%%%%%%%%%%%%%%%%%%%%%%%%%%%%%%%%%%%%%%%%%%%%%%%%%
\section{Introduction}

\noindent Observability and controllability are two fundamental structural properties of a control system.
The former describes the possibility of inferring the state that characterizes the system from observing its inputs and outputs. The latter characterizes the possibility to move a system in all its space of states, by using suitable system inputs (controls).
Both these concepts were first introduced for linear systems \cite{Kalman61,Kalman63} and
the two analytic conditions to check if a linear system satisfies these two properties have also be obtained.
\noindent The nonlinear case is much more complex. First, both these concepts become local. In addition, unlike the linear case, observability depends on the system inputs.
In this paper we refer to the {\it weak local observability}, as defined in \cite{Her77,Casti82}
(definitions 8, 9, 10, 11, in \cite{Casti82}).
Regarding controllability, in this paper we refer to the concept of {\it weak local controllability}, as defined in \cite{Her77}.

\noindent The two analytic conditions to check if a continuous time-invariant nonlinear system satisfies these two properties (the weak local observability and the weak local controllability) have also been introduced \cite{Her77,Casti82,Stefan74,Suss83,Isi95,Lewis01}. They are known as the {\it observability rank condition} and the {\it controllability rank condition}. They are summarized in section \ref{SubSectionObsTI} and in section \ref{SubSectionConTI}, respectively.
Very recently, new analytic conditions have also been proposed.
The conditions proposed in \cite{TAC19,ARCH17b}, extend the observability rank condition to the case when the dynamics are also driven by unknown inputs. The authors of \cite{Carravetta19} proposed a new condition for the weak local controllability, which presents some advantages with respect to the controllability rank condition.

\noindent Unfortunately, all the conditions above cannot be used in the case when the system is time-varying (nonautonomous).

\noindent A time-varying system is a system whose behaviour changes with time. In particular, the system will respond differently to the same input at different times.
{\blue A typical example of time-varying system is an aircraft. For this system, there are two main factors that make it time-varying: decreasing weight due to consumption of fuel and the different configuration of control surfaces during take off, cruise and landing. The first factor will characterize the system investigated in sections \ref{SectionApplicationAerospace} and \ref{SectionTakeOff}.}

\noindent In a general mathematical characterization of a time-varying nonlinear system, all the key scalar and vector fields that define its dynamics and/or its output functions, explicitly depend upon time (see equation (\ref{EquationObsSystemDefinitionTV})).

\noindent For time-varying systems, the two analytic conditions to check observability and controllability have only been obtained in the linear case. These conditions are summarized in section \ref{SubSectionObsTVLinear} and in section \ref{SubSectionConTVLinear}, respectively. 

\noindent So far, no condition exists to check 
the weak local observability and the weak local controllability for time-varying nonlinear systems. This is precisely the goal of this paper.

\noindent Specifically, the contributions of this paper are the following two:

\begin{enumerate}

\item Extend the observability rank condition to nonlinear time-varying systems.

\item Extend the controllability rank condition to nonlinear time-varying systems.

\end{enumerate}

%%%RRR\noindent {\blue The two new conditions consist in the calculation of a codistribution and of a distribution (respectively for observability and controllability), rather than in the calculation of the rank of matrices, as in the linear case.}

\noindent The paper is organized as follows. Section \ref{SectionSystem} provides the basic equations that characterize the systems here investigated. Sections \ref{SectionObservability} and \ref{SectionControllability} provide the two new analytic conditions, whose derivations are given separately in section \ref{SectionDerivation}.
Section \ref{SectionExamples} provides two simple applications. {\blue They are deliberately trivial to better illustrate the two new analytic conditions. Sections  \ref{SectionApplicationAerospace} and 
\ref{SectionTakeOff}
provide a real application. We investigate the observability and controllability properties of a lunar module that operates in presence of gravity and in absence of an atmosphere. This system has an explicit time dependence due to the fuel consumption that results in a variation of the weight and the variation of the moment of inertia.} Finally, our conclusion is given in section \ref{SectionConclusion}.

%\ref{AlgoObsTI} 1
%\ref{AlgoObsTV} 2
%\ref{AlgoConTI} 3
%\ref{AlgoConTV} 4
%\ref{AlgoObsTVC} 5
%\ref{AlgoConTVC} 6

\section{Considered systems}\label{SectionSystem}

\noindent We will refer to a nonlinear control system with $m$ inputs ($ u_1,\cdots,u_m$). The state is the vector

\[
x\triangleq [x^1,\cdots,~x^n]^T \in \mathcal{M}
\]

\noindent with $\mathcal{M}$ an open set of $\mathbb{R}^n$. We assume that the dynamics are nonlinear with respect to the state and affine with respect to the inputs. 
We account for an explicit time dependence, namely, all the functions that characterize the dynamics and/or the outputs, can explicitly depend on time. 
Finally, the system has $p(\ge 1)$ outputs. 
Our system is characterized by the following equations:

\begin{equation}\label{EquationObsSystemDefinitionTV}
\left\{\begin{array}{ll}
  \dot{x} &=   f^0(x,~t)  + \sum_{i=1}^m f^i (x,~t) u_i\\
  y &= [h_1( x,~t),\cdots,h_p(x,~t)]^T \\
\end{array}\right.
\end{equation}

\noindent where $ f^i (x,~t)$, $i=0,1,\cdots,m$, are vector fields in $\mathcal{M}$  and the functions $h_1(x,~t), ~\cdots, ~h_p(x,~t)$ are scalar fields defined on the open set $\mathcal{M}$. All these vector and scalar fields explicitly depend on time.

\section{Analytic condition for observability}\label{SectionObservability}

\noindent This section introduces the analytic condition to check the state observability for systems that satisfy equation (\ref{EquationObsSystemDefinitionTV}).

\noindent Before introducing this new condition we remind the reader the existing results for the less general systems. Specifically, in section \ref{SubSectionObsTVLinear}, we provide the analytic condition that holds in the case of time-varying linear systems and, in section \ref{SubSectionObsTI}, we provide the analytic condition that holds in the case of time-invariant nonlinear systems. In section 
\ref{SubSectionObsTV}, we provide the new condition that holds in general, i.e., for  time-varying nonlinear systems.

\subsection{Time-varying linear systems}\label{SubSectionObsTVLinear}

\noindent This special case is obtained by setting in (\ref{EquationObsSystemDefinitionTV}):

\begin{itemize}

\item $f^0(x,~t)=A(t)x$, where $A$ is a matrix of dimension $n\times n$.

\item $f^i(x,~t)=b^i(t)$, where $b^1(t),\cdots,b^m(t)$ are $m$ column vectors of dimension $n$.

\item $h_j(x,~t)=c_j(t)x$, where $c_1(t),\cdots,c_p(t)$ are $p$ row-vectors of dimension $n$.

\end{itemize}

\noindent We can write (\ref{EquationObsSystemDefinitionTV}) as follows:

\begin{equation}\label{EquationObsSystemDefinitionTVLinear}
\left\{\begin{array}{ll}
  \dot{x} &=   A(t)x  + B(t) u\\
  y &= C(t)x \\
\end{array}\right.
\end{equation}

\noindent where the columns of $B$ are the vectors $b^1,\cdots,b^m$ above and the lines of $C$ are the vectors $c_1,\cdots,c_p$ above.

\noindent The system defined by (\ref{EquationObsSystemDefinitionTVLinear}) is observable in a given time interval $\mathcal{I}$ if there exists $\bar{t}\in\mathcal{I}$ and a positive integer $k$ such that:

\begin{equation}\label{EquationObsCondLin}
rank\left[
\begin{array}{c}
  N_0(\bar{t}) \\
  N_1(\bar{t}) \\
  \cdots \\
  N_k(\bar{t}) \\
\end{array}
\right]=n
\end{equation}

\noindent where $N_0(t)\triangleq C(t)$ and $N_i(t)$  is defined recursively as:

\begin{equation}\label{EquationObsCodLin}
N_i(t)=N_{i-1}(t)A(t) + \frac{dN_{i-1}(t)}{dt}, ~~i=1,\cdots,k
\end{equation} 

{\blue
\noindent This result was obtained long time ago in \cite{Sil67}. }The reader is also addressed to 
\cite{Sontag13} for further details and for the analytic derivations to prove the validity of the above condition.

\subsection{Time-invariant nonlinear systems}\label{SubSectionObsTI}

\noindent This special case is obtained when all the vector and scalar fields that appear in (\ref{EquationObsSystemDefinitionTV}) do not explicitly depend on time.
The analytic condition to check the weak local observability at a given $x_0 \in \mathcal{M}$ of the state $x$ that satisfies (\ref{EquationObsSystemDefinitionTV})
is obtained by computing the {\it observable codistribution} \cite{Isi95}. When all the vector and scalar fields do not explicitly depend on time, the observable codistribution is generated by the recursive algorithm 1 (see \cite{Her77,Casti82,Isi95}). We use the following notation:

\begin{itemize}

\item Given a scalar field $h$, $dh$ is its differential.

\item Given a vector field $f$ (defined on the open set $\mathcal{M}$),  $\mathcal{L}_f$ denotes the Lie derivative along $f$. We remind the reader that, the Lie derivative along $f$ of a given scalar field $h$ is \cite{Isi95}:

\[
\mathcal{L}_fh=\frac{\partial h}{\partial x} \cdot f
\]
\noindent Additionally:

\begin{equation}
\label{EquationLieDifferential}
\mathcal{L}_fdh=d\mathcal{L}_fh
\end{equation}

\item Given a codistribution $\Omega$ and a given vector field $f$ (both defined on the open set $\mathcal{M}$), $\mathcal{L}_f \Omega$ denotes the codistribution whose covectors are the Lie derivatives along $f$ of the covectors in $\Omega$.

\item Given two vector spaces $V_1$ and $V_2$,  $V_1{\blue\oplus}V_2$ is their sum, i.e., the span of all the generators of $V_1$ and $V_2$. 

\end{itemize}

\begin{algorithm}
\begin{algorithmic}\label{AlgoObsTI}
{\blue
  \STATE Set  $\Omega=$span$\left\{d  h_1,\cdots, d h_p \right\}$ 
  \WHILE{dim$\left(\Omega\oplus\mathcal{L}_{f^0}\Omega\oplus\cdots\oplus\mathcal{L}_{f^m}\Omega\right)>$dim$(\Omega)$}
   \STATE Set 
   $\Omega=\Omega\oplus\mathcal{L}_{f^0}\Omega\oplus\cdots\oplus\mathcal{L}_{f^m}\Omega$
  \ENDWHILE}
\end{algorithmic}
\caption{Observable codistribution for time-invariant nonlinear systems.}
\end{algorithm}

%\begin{algorithm}[h]
%\begin{algorithmic}\label{AlgoObsTI}
%  \State Set $k=0$
%  \State Set $\Omega_k=span\left\{d  h_1,\cdots, d h_p \right\}$ 
%  \State Set $k=k+1$
%  \State Set $\Omega_k=\Omega_{k-1}+\sum_{i=0}^m\mathcal{L}_{f^i}\Omega_{k-1}$
%  \While{$dim(\Omega_k)>dim(\Omega_{k-1})$}
%   \State Set $k=k+1$
%  \State Set $\Omega_k=\Omega_{k-1}+\sum_{i=0}^m\mathcal{L}_{f^i}\Omega_{k-1}$
%  \EndWhile
%  \State Set $\Omega=\Omega_{k-1}$ and $s=dim(\Omega)$
%\end{algorithmic}
%\caption{Observable codistribution for time-invarying nonlinear systems.}
%\end{algorithm}

\noindent The analytic condition to check the weak local observability of nonlinear time-invariant systems is given by the following fundamental result:

\begin{theorem}[Observability Rank Condition]\label{TheoremORC}
Algorithm 1 converges in an open and dense set of $\mathcal{M}$ and the convergent codistribution is obtained {\blue in at most $ n-1$ steps}.
If the convergent codistribution is non singular at $x_0\in\mathcal{M}$ and its dimension is equal to $n$ at $x_0$, then the system is weakly locally observable at $x_0$ (sufficient condition). Conversely, 
if the system is weakly locally observable at $x_0$, the dimension of the above codistribution is $n$ in a dense neighbourhood of $x_0$ 
(necessary condition).
\end{theorem}

\proof{All the statements are very well known results. The reader is addressed to \cite{Isi95} (lemmas 1.9.1, 1.9.2 and 1.9.6) for the convergence properties of algorithm 1. The proof of the sufficient condition is available in \cite{Her77}, theorem 3.1. The proof of the necessary condition is available in  \cite{Her77}, theorem 3.11
$\blacktriangleleft$}

\subsection{Time-varying nonlinear systems}\label{SubSectionObsTV}

\noindent We now consider the general case of time-varying nonlinear systems. In this section we only provide the analytic condition. In section \ref{SubSectionDerivationObs} we prove its validity.

\noindent The new condition is similar to the condition that holds in the case time-invariant (i.e., the observability rank condition provided in section \ref{SubSectionObsTI}). The only difference resides in the computation of the observable codistribution.
The new codistribution is given by algorithm 2,
where we introduced the following operator:

\begin{equation}\label{EquationMixOperatorTV}
\widetilde{\mathcal{L}}_{f^0} \triangleq \frac{\partial}{\partial t} + \mathcal{L}_{f^0}
\end{equation}

\begin{algorithm}[h]
\begin{algorithmic}\label{AlgoObsTV}
{\blue
  \STATE Set  $\Omega=$span$\left\{d  h_1,\cdots, d h_p \right\}$ 
  \WHILE{dim$\left(\Omega\oplus\widetilde{\mathcal{L}}_{f^0}\Omega\oplus\mathcal{L}_{f^1}\Omega\oplus\cdots\oplus\mathcal{L}_{f^m}\Omega\right)>$dim$(\Omega)$}
   \STATE Set 
   $\Omega=\Omega\oplus\widetilde{\mathcal{L}}_{f^0}\Omega\oplus\mathcal{L}_{f^1}\Omega\oplus\cdots\oplus\mathcal{L}_{f^m}\Omega$
  \ENDWHILE}
\end{algorithmic}
\caption{Observable codistribution for time-variant nonlinear systems.}
\end{algorithm}

%\begin{algorithm}
%\begin{algorithmic}\label{AlgoObsTV}
%  \STATE Set $k=0$
%  \STATE Set $\Omega_k=span\left\{d  h_1,\cdots, d h_p \right\}$ 
%  \STATE Set $k=k+1$
%  \STATE Set $\Omega_k=\Omega_{k-1}+\widetilde{\mathcal{L}}_{f^0}\Omega_{k-1}+\sum_{i=1}^m\mathcal{L}_{f^i}\Omega_{k-1}$
%  \WHILE{$dim(\Omega_k)>dim(\Omega_{k-1})$}
%   \STATE Set $k=k+1$
%  \STATE Set $\Omega_k=\Omega_{k-1}+\widetilde{\mathcal{L}}_{f^0}\Omega_{k-1}+\sum_{i=1}^m\mathcal{L}_{f^i}\Omega_{k-1}$
%  \ENDWHILE
%  \STATE Set $\Omega=\Omega_{k-1}$ and $s=dim(\Omega)$
%\end{algorithmic}
%\caption{Observable codistribution for time-varying nonlinear systems.}
%\end{algorithm}
%

%\begin{algorithm}[h]
%\begin{algorithmic}\label{AlgoObsTV}
%  \State Set $k=0$
%  \State Set $\Omega_k=span\left\{d  h_1,\cdots, d h_p \right\}$ 
%  \State Set $k=k+1$
%  \State Set $\Omega_k=\Omega_{k-1}+\widetilde{\mathcal{L}}_{f^0}\Omega_{k-1}+\sum_{i=1}^m\mathcal{L}_{f^i}\Omega_{k-1}$
%  \While{$dim(\Omega_k)>dim(\Omega_{k-1})$}
%   \State Set $k=k+1$
%  \State Set $\Omega_k=\Omega_{k-1}+\widetilde{\mathcal{L}}_{f^0}\Omega_{k-1}+\sum_{i=1}^m\mathcal{L}_{f^i}\Omega_{k-1}$
%  \EndWhile
%  \State Set $\Omega=\Omega_{k-1}$ and $s=dim(\Omega)$
%\end{algorithmic}
%\caption{Observable codistribution for time-varying nonlinear systems.}
%\end{algorithm}

%\noindent The derivations provided in section \ref{SectionDerivation} show that the convergence properties of algorithm 2 are the same of algorithm 1. In particular, 
%the convergent codistribution is obtained at the smallest integer $m$ for which $\Omega_m=\Omega_{m-1}$ and $m \le n-1$.

\noindent Note that the codistribution returned by the algorithm above is in general time-dependent. 

\noindent The analytic condition to check the weak local observability of nonlinear time-varying systems is given by the following fundamental new result:

\begin{theorem}[Extended Observability Rank Condition]\label{TheoremEORC}
Algorithm 2 converges in an open and dense set of {\blue $\mathbb{R}\times\mathcal{M}$} and the convergent codistribution is obtained {\blue in at most $n-1$ steps}.
If the convergent codistribution is non singular at $x_0\in\mathcal{M}$ and at a given time $t_0\in\mathbb{R}$ and its dimension is equal to $n$ at $(t_0,~x_0)$, then the system is weakly locally observable at $(t_0,~x_0)$ (sufficient condition). Conversely, 
if the system is weakly locally observable at $(t_0,~x_0)$, the dimension of the above codistribution is $n$ in a dense neighbourhood of $(t_0,~x_0)$
(necessary condition).
\end{theorem}

\proof{The proof is given in section \ref{SubSectionDerivationObs}
$\blacktriangleleft$}

\vskip.3cm

\noindent We conclude this section with the following remarks:

\begin{enumerate}

\item Algorithm 2 differs from algorithm 1 only for the recursive step. In particular, the operator given in (\ref{EquationMixOperatorTV}) substitutes the Lie derivative along $f^0$. In other words, the new algorithm is obtained with the substitution:

\[
\mathcal{L}_{f^0}\rightarrow
\widetilde{\mathcal{L}}_{f^0}
\]

{\blue

\noindent If $f^0$ is null, in the recursive step we need to add the term $\oplus \frac{\partial}{\partial t}\Omega$.

\item If $\Omega$ is generated by $\omega_1,\cdots,\omega_s$, the codistribution $\Omega\oplus
\widetilde{\mathcal{L}}_{f^0}\Omega\oplus\cdots\oplus\mathcal{L}_{f^m}\Omega$
is generated by $\omega_1,..,\omega_s,\widetilde{\mathcal{L}}_{f^0}\omega_1,..,\widetilde{\mathcal{L}}_{f^0}\omega_s,..,\mathcal{L}_{f^m}\omega_1,..,\mathcal{L}_{f^m}\omega_s$ (see appendix \ref{AppendixCod}). This allows us to easily implement algorithm 2 since it suffices to compute the Lie derivatives of the generators of $\Omega$, at each step.
}

\item The extended observability rank condition reduces to the observability rank condition provided in section \ref{SubSectionObsTI} when all the vector and scalar fields that appear in (\ref{EquationObsSystemDefinitionTV}) do not explicitly depend on time.

\item The extended observability rank condition reduces to the condition provided in section \ref{SubSectionObsTVLinear} in the linear case.

\end{enumerate}

\section{Analytic condition for controllability}\label{SectionControllability}

\noindent This section introduces the analytic condition to check the state controllability for systems that satisfy equation (\ref{EquationObsSystemDefinitionTV}). The section is structured as the section \ref{SectionObservability}.
Specifically, in section \ref{SubSectionConTVLinear}, we provide the analytic condition that holds in the case of time-varying linear systems and, in section \ref{SubSectionConTI}, we provide the analytic condition that holds in the case of time-invariant nonlinear systems. Finally, in section 
\ref{SubSectionConTV}, we provide the new condition that holds in general, i.e., for  time-varying nonlinear systems.

\subsection{Time-varying linear systems}\label{SubSectionConTVLinear}

\noindent This special case is the same considered in section \ref{SubSectionObsTVLinear}. In other words, we refer to the system characterized by (\ref{EquationObsSystemDefinitionTVLinear}).
This system is controllable if there exists $\bar{t}\in\mathcal{I}$ and a positive integer $k$ such that:

\begin{equation}\label{EquationConCondLin}
rank\left[
M_0, ~M_1, \cdots, ~M_k
\right]=n
\end{equation}

\noindent where $M_0(t)\triangleq B(t)$ and $M_i(t)$  is defined recursively as:

\begin{equation}\label{EquationConDistLin}
M_i(t)=A(t)M_{i-1}(t) - \frac{dM_{i-1}(t)}{dt}, ~~i=1,\cdots,k
\end{equation} 

{\blue
\noindent This result was obtained long time ago in \cite{Sil67}. }The reader is also addressed to 
\cite{Sontag13} for further details and for the analytic derivations to prove the validity of the above condition.

\subsection{Time-invariant nonlinear systems}\label{SubSectionConTI}

\noindent This special case is obtained when all the vector and scalar fields that appear in (\ref{EquationObsSystemDefinitionTV}) do not explicitly depend on time.
The analytic condition to check the weak local controllability from a given $x_0 \in \mathcal{M}$
is obtained by computing the {\it controllability distribution} \cite{Isi95}. 
When all the vector and scalar fields do not explicitly depend on time, the controllability distribution is generated by the recursive algorithm 3 (see \cite{Isi95}).
We use the following notation:

\begin{itemize}

\item Given two vector fields $f^a, ~f^b$ (defined on the open set $\mathcal{M}$),  $[f^a,~f^b]$ denotes their Lie bracket, defined as follows:

\begin{equation}
\label{EquationLieBracket}
[f^a,~f^b]=\frac{\partial f^b}{\partial x}f^a-\frac{\partial f^a}{\partial x}f^b
\end{equation}

\item Given a distribution $\Delta$ and a given vector field $f$ (both defined on the open set $\mathcal{M}$), $[\Delta, ~f]$ denotes the distribution whose vectors are the Lie bracket of any vector in $\Delta$ with $f$.

\end{itemize}

\begin{algorithm}[h]
\begin{algorithmic}\label{AlgoConTI}
{\blue
  \STATE Set $\Delta=$span$\left\{f^1,\cdots, f^m \right\}$ 
  \WHILE{dim$\left(\Delta\oplus [\Delta,~f^0]\oplus\cdots\oplus [\Delta,~f^m]\right)>$dim$(\Delta)$}
  \STATE Set $\Delta=\Delta\oplus [\Delta,~f^0]\oplus\cdots\oplus [\Delta,~f^m]$
  \ENDWHILE}
\end{algorithmic}
\caption{Controllability distribution for time-invariant nonlinear systems.}
\end{algorithm}

%\begin{algorithm}[h]
%\begin{algorithmic}\label{AlgoConTI}
%  \State Set $k=0$
%  \State Set $\Delta_k=span\left\{f^1,\cdots, f^m \right\}$ 
%  \State Set $k=k+1$
%  \State Set $\Delta_k=\Delta_{k-1}+\sum_{i=0}^m [\Delta_{k-1},f^i]$
%  \While{$dim(\Delta_k)>dim(\Delta_{k-1})$}
%   \State Set $k=k+1$
%  \State Set $\Delta_k=\Delta_{k-1}+\sum_{i=0}^m [\Delta_{k-1},f^i]$
%  \EndWhile
%  \State Set $\Delta=\Delta_{k-1}$ and $d=dim(\Delta)$
%\end{algorithmic}
%\caption{Controllability distribution for time-invariant nonlinear systems.}
%\end{algorithm}
%

\noindent The analytic condition to check the weak local controllability of nonlinear time-invariant systems is given by the following fundamental result:

\begin{theorem}[Controllability Rank Condition]\label{TheoremCRC}
Algorithm 3 converges in an open and dense set of $\mathcal{M}$ and the convergent distribution is obtained {\blue in at most $n-1$ steps}.
If the convergent distribution is non singular at $x_0\in\mathcal{M}$ and its dimension is equal to $n$ at $x_0$, then the system is 
weakly locally controllable from $x_0$ (sufficient condition). Conversely, 
if the system is weakly locally controllable from $x_0$, the dimension of the above distribution is $n$ in a dense neighbourhood of $x_0$ (necessary condition).
\end{theorem}

\proof{All the statements are very well known results. The reader is addressed to \cite{Isi95} (lemmas 1.8.1, 1.8.2 and 1.8.3) for the convergence properties of algorithm 3. The proof of the sufficient condition is available in \cite{Her77}, theorem 2.2. The proof of the necessary condition is available in  \cite{Her77}, theorem 2.5
$\blacktriangleleft$}

\subsection{Time-varying nonlinear systems}\label{SubSectionConTV}

\noindent We now consider the general case of time-varying nonlinear systems. In this section we only provide the analytic condition. In section \ref{SubSectionDerivationCon} we prove its validity.

\noindent The new condition is similar to the condition that holds in the case time-invariant (i.e., the controllability rank condition provided in section \ref{SubSectionConTI}). The only difference resides in the computation of the controllability distribution.

\noindent The new distribution is given by algorithm 4,
where we introduced the following operator:

\begin{equation}\label{EquationMixBracketTV}
\langle a, ~f^0\rangle \triangleq [a, ~f^0] - \frac{\partial a}{\partial t}
\end{equation}

\begin{algorithm}[h]
\begin{algorithmic}\label{AlgoConTV}
{\blue
  \STATE Set $\Delta=$span$\left\{f^1,\cdots, f^m \right\}$ 
  \WHILE{dim$\left(\Delta\oplus\langle \Delta, f^0\rangle\oplus [\Delta,f^1]\oplus ..\oplus [\Delta,f^m]\right)>$dim$(\Delta)$}
  \STATE Set $\Delta=\Delta\oplus\langle \Delta, ~f^0\rangle\oplus [\Delta,~f^1]\oplus\cdots\oplus [\Delta,~f^m]$
  \ENDWHILE}
\end{algorithmic}
\caption{Controllability distribution for time-variant nonlinear systems.}
\end{algorithm}

%\begin{algorithm}[h]
%\begin{algorithmic}\label{AlgoConTV}
%  \State Set $k=0$
%  \State Set $\Delta_k=span\left\{f^1,\cdots, f^m \right\}$ 
%  \State Set $k=k+1$
%  \State Set $\Delta_k=\Delta_{k-1}+\sum_{i=1}^m [\Delta_{k-1},f^i]+\langle\Delta, ~f^0\rangle$
%  \While{$dim(\Delta_k)>dim(\Delta_{k-1})$}
%   \State Set $k=k+1$
%  \State Set $\Delta_k=\Delta_{k-1}+\sum_{i=1}^m [\Delta_{k-1},f^i]+\langle\Delta, ~f^0\rangle$
%  \EndWhile
%  \State Set $\Delta=\Delta_{k-1}$ and $d=dim(\Delta)$
%\end{algorithmic}
%\caption{Controllability distribution for time-varying nonlinear systems.}
%\end{algorithm}

%\noindent The derivations provided in section \ref{SectionDerivation} show that the convergence properties of algorithm 4 are the same of algorithm 3. In particular, 
%the convergent distribution is obtained at the smallest integer $k$ for which $\Delta_k=\Delta_{k-1}$ and $k \le n-1$.

\noindent Note that the distribution returned by the algorithm above is in general time-dependent. 

\noindent The analytic condition to check the weak local controllability of nonlinear time-varying systems is given by the following fundamental new result:

\begin{theorem}[Extended Controllability Rank Condition]\label{TheoremECRC}
Algorithm 4 converges in an open and dense set of {\blue $\mathbb{R}\times\mathcal{M}$} and the convergent distribution is obtained {\blue in at most $n-1$ steps}.
If the convergent distribution is non singular at $x_0\in\mathcal{M}$ and at a given time $t_0\in\mathbb{R}$ and its dimension is equal to $n$ at $(t_0,~x_0)$, then the system is weakly locally controllable from $(t_0,~x_0)$ (sufficient condition). Conversely, 
if the system is weakly locally controllable from $(t_0,~x_0)$, the dimension of the above distribution is $n$ in a dense neighbourhood of $(t_0,~x_0)$
(necessary condition).
\end{theorem}

\proof{The proof is given in section \ref{SubSectionDerivationCon}
$\blacktriangleleft$}

\vskip.3cm

\noindent We conclude this section with the following remarks:

\begin{enumerate}

\item Algorithm 4 differs from algorithm 3 only for the recursive step. In particular, the operator given in (\ref{EquationMixBracketTV}) substitutes the Lie bracket with $f^0$. In other words, the new algorithm is obtained with the substitution:

\[
[\cdot, ~f^0]\rightarrow\langle\cdot, ~f^0\rangle
\]

{\blue

\noindent If $f^0$ is null, in the recursive step we need to add the term $\oplus \frac{\partial}{\partial t}\Delta$.

\item If $\Delta$ is generated by $d^1,\cdots,d^k$, the distribution $\Delta\oplus
\langle\Delta,~f^0\rangle\oplus\cdots\oplus
\langle\Delta,~f^m\rangle$
is generated by $d^1,..,d^k,\langle d^1,~f^0\rangle,..,\langle d^k,~f^0\rangle,.. [d^1,~f^m],..,[d^k,~f^m]$  (see appendix \ref{AppendixDist}). This allows us to easily implement algorithm 4 since it suffices to compute the Lie brackets of the generators of $\Delta$, at each step.
}

\item The extended controllability rank condition reduces to the controllability rank condition provided in section \ref{SubSectionConTI} when all the vector and scalar fields that appear in (\ref{EquationObsSystemDefinitionTV}) do not explicitly depend on time.

\item The extended controllability rank condition reduces to the condition provided in section \ref{SubSectionConTVLinear} in the linear case.

\end{enumerate}

\section{Proofs}\label{SectionDerivation}

\noindent In this section we prove the validity of the theorems \ref{TheoremEORC} and \ref{TheoremECRC} (sections \ref{SubSectionDerivationObs} and \ref{SubSectionDerivationCon}, respectively).

\noindent Both these proofs are obtained by including in the state the variable time.
We denote the new extended state by $\underline{x}$ and we have
$\underline{x}=\left[t,~x^1,~\cdots,~x^n\right]^T$.

\noindent We characterize the system in (\ref{EquationObsSystemDefinitionTV}) by using the extended state. From the first equation in (\ref{EquationObsSystemDefinitionTV}) we obtain:

\begin{equation}\label{EquationChronoDynamics}
  \dot{\underline{x}} =   \underline{f}^0(\underline{x}) +  \sum_{i=1}^m \underline{f}^i(\underline{x})u_i
\end{equation}

%\noindent with:

\begin{equation}\label{EquationG0F}
\underline{f}^0(\underline{x}) \equiv
\left[
\begin{array}{c}
1 \\
f^0
\end{array}
\right],~~~~
\underline{f}^i(\underline{x}) \equiv
\left[
\begin{array}{c}
0 \\
f^i
\end{array}
\right]
\end{equation}

\noindent Regarding the outputs, we remark that we need to include a new output that is $h_0(\underline{x})=t$. Indeed, it is a common (and implicit) assumption that all the system inputs and the outputs are synchronized. For instance, in a real system, the inputs and outputs are measured by sensors. The sensors provide their measurements together with the time when each measurement has occurred.
This means that our system is also equipped with an additional sensor that is the clock (i.e., a sensor that measures time). Therefore, a full description of our system in the extended state is given by:

\begin{equation}\label{EquationChronoDynamicsOutputsTV}
\left\{\begin{array}{ll}
  \dot{\underline{x}} &=   \underline{f}^0(\underline{x}) +  \sum_{i=1}^m \underline{f}^i(\underline{x})u_i\\
  y &= [h_0(\underline{x})=t, ~h_1(\underline{x}),\cdots,h_p(\underline{x})]^T \\
\end{array}\right.
\end{equation}

%
%\noindent with $h_0(\underline{x})=t$.

\subsection{Proof of theorem \ref{TheoremEORC}}\label{SubSectionDerivationObs}

%\proof{

\noindent By introducing the extended state we transformed our original nonautonomous system in (\ref{EquationObsSystemDefinitionTV}) into the autonomous system in (\ref{EquationChronoDynamicsOutputsTV}). We are allowed to use the results stated by theorem \ref{TheoremORC}.
The algorithm that provides the 
observable codistribution in the extended state is algorithm 5.

\begin{algorithm}[h]
\begin{algorithmic}
{\blue
  \STATE Set $\underline{\Omega}=$span$\left\{\underline{d}h_0, \underline{d}h_1,\cdots,\underline{d}h_p\right\}$ 
  \WHILE{dim$\left(\underline{\Omega}\oplus\mathcal{L}_{\underline{f}^0}\underline{\Omega}\oplus\cdots\oplus\mathcal{L}_{\underline{f}^m}\underline{\Omega}\right)>$dim$(\underline{\Omega})$}
  \STATE Set $\underline{\Omega}=\underline{\Omega}\oplus\mathcal{L}_{\underline{f}^0}\underline{\Omega}\oplus\cdots\oplus\mathcal{L}_{\underline{f}^m}\underline{\Omega}$
  \ENDWHILE}
\end{algorithmic}\label{AlgoObsTVC}
\caption{Observable codistribution in the extended state}
\end{algorithm}

%\begin{algorithm}[h]
%\begin{algorithmic}
%  \State Set $k=0$
%  \State Set $\underline{\Omega}_k=span\left\{\underline{d}h_0, \underline{d}h_1,\cdots,\underline{d}h_p\right\}$ 
%  \State Set $k=k+1$
%  \State Set $\underline{\Omega}_k=\underline{\Omega}_{k-1}+\mathcal{L}_{\underline{f}^0}\underline{\Omega}_{k-1}+\sum_{i=1}^m\mathcal{L}_{\underline{f}^i}\underline{\Omega}_{k-1}$
%  \While{$dim(\underline{\Omega}_k)>dim(\underline{\Omega}_{k-1})$}
%   \State Set $k=k+1$
%  \State Set $\underline{\Omega}_k=\underline{\Omega}_{k-1}+\mathcal{L}_{\underline{f}^0}\underline{\Omega}_{k-1}+\sum_{i=1}^m\mathcal{L}_{\underline{f}^i}\underline{\Omega}_{k-1}$
%  \EndWhile
%\end{algorithmic}\label{AlgoObsTVC}
%\caption{Observable codistribution in the extended state}
%\end{algorithm}
%
%

\noindent We denoted by $\underline{d}$ the differential in the extended state. From theorem \ref{TheoremORC} we know that algorithm 5  converges 
in an open and dense set of {\blue $\mathbb{R}\times\mathcal{M}$} and the convergent codistribution is obtained {\blue in at most $(n+1)-1(=n)$ steps}.
If the convergent codistribution is non singular at $(t_0,~x_0)$ and its dimension is equal to $n+1$ at $(t_0,~x_0)$, then the system is weakly locally observable at $(t_0,~x_0)$. Conversely, 
if the system is weakly locally observable at $(t_0,~x_0)$, the dimension of the above codistribution is $n+1$ in a dense neighbourhood of $(t_0,~x_0)$.

\noindent From these results, we immediately obtain the proof of the results stated by theorem \ref{TheoremEORC} by using the following 
fundamental separation property. {\blue At every step of algorithm 5, we can 
split $\underline{\Omega}$} into two codistributions, as follows:

\begin{equation}\label{EquationSeparation}
\underline{\Omega}= \textnormal{span}\{dt\} \oplus \Omega
\end{equation}

\noindent where $\Omega$ is generated by differentials of scalar fields only with respect to the state (and not the extended state). The validity of the property above is a consequence of the fact that 
the extended system is characterized by the output $h_0=t$ and, consequently, we have $\underline{d}h_0=dt$.

%\noindent Hence, $\Omega_k$ is precisely the observable codistribution we want to compute, i.e., the one that only includes the observability properties of the state. Our objective is to derive the algorithm that directly provides the above $\Omega_k$.
%We proceed as follows.

\noindent For any $h$ such that $\underline{d}h\in \underline{\Omega}$ {\blue at a given step}, the following $m+1$ covectors belong to $\underline{\Omega}$ {\blue at the next step}:

%\[
%\begin{array}{l}
%\mathcal{L}_{\underline{f}^0}~\underline{d}h  \\
% \mathcal{L}_{\underline{f}^1}~\underline{d}h\\
% \cdots\\
% \mathcal{L}_{\underline{f}^m}~\underline{d}h\\
%\end{array}
%\]
%

\[
\mathcal{L}_{\underline{f}^0}~\underline{d}h,~~
 \mathcal{L}_{\underline{f}^1}~\underline{d}h,~~
 \cdots,~~
 \mathcal{L}_{\underline{f}^m}~\underline{d}h
\]

\noindent From the structure of $\underline{f}^j$ given in (\ref{EquationG0F}) we obtain:

\begin{equation}\label{EquationLieDevelopedTVC}
\mathcal{L}_{\underline{f}^j}h= 
\left[\begin{array}{ll}
  \frac{\partial h}{\partial t} + \mathcal{L}_{f^0}h &j=0\\
  \mathcal{L}_{f^j}h &j=1,\cdots,m\\
\end{array}
\right.
\end{equation}

\noindent As result, by using (\ref{EquationLieDifferential}), for any $h$ such that $\underline{d}h\in \underline{\Omega}$ {\blue at a given step}, the following $m+1$ covectors belong to $\underline{\Omega}$ {\blue at the next step}:

%\[
%\begin{array}{l}
%\underline{d}\widetilde{\mathcal{L}}_{f^0}h =\underline{d}\left(\frac{\partial h}{\partial t} + \mathcal{L}_{f^0}h\right) \\
% \underline{d}\mathcal{L}_{f^1}h\\
% \cdots\\
% \underline{d}\mathcal{L}_{f^m}h\\
%\end{array}
%\]
%

\[
\underline{d}\widetilde{\mathcal{L}}_{f^0}h =\underline{d}\left(\frac{\partial h}{\partial t} + \mathcal{L}_{f^0}h\right),~~
 \underline{d}\mathcal{L}_{f^1}h,~~
 \cdots,~~
 \underline{d}\mathcal{L}_{f^m}h
\]

\noindent Finally, by using (\ref{EquationSeparation}), we obtain that, 
for any $h$ such that $dh\in\Omega$ {\blue at a given step}, the following $m+1$ covectors belong to $\Omega$ {\blue at the next step}:

%\[
%\begin{array}{l}
%d\left(\frac{\partial h}{\partial t} + \mathcal{L}_{f^0}h\right) \\
%d\mathcal{L}_{f^1}h\\
% \cdots\\
%d\mathcal{L}_{f^m}h\\
%\end{array}
%\]

\[
d\left(\frac{\partial h}{\partial t} + \mathcal{L}_{f^0}h\right),~~
d\mathcal{L}_{f^1}h,~~
 \cdots,~~
d\mathcal{L}_{f^m}h
\]

\noindent This proves that $\Omega$ is generated by algorithm 2.

\noindent Finally,  the convergence of algorithm 5 (and consequently of algorithm 2) occurs {\blue in at most $n-1$ steps instead of $(n+1)-1(=n)$ steps}.
This is proved as follows. 
The dimension of $\underline{\Omega}$ {\blue at the initialization} satisfies:
\[
\dim\left(\underline{\Omega}\right)\ge2
\]
\noindent in an open and dense set of $\mathbb{R}\times\mathcal{M}$. Indeed, at the initialization,
$\dim\left(\underline{\Omega}\right)=\dim($span$\{dt\} + \Omega )= $
$\dim\left(\right.$span$\left.\{dt\} +\right.$span$\left.\{dh_1,~\cdots,dh_p\} \right)=$
$\dim\left(\right.$span$\left.\{dt\} \right) +\dim\left(\right.$span$\left.\{dh_1,~\cdots,dh_p\} \right)=$
$
1+\dim\left(\right.$span$\left.\{h_1,~\cdots,h_p\}\right)\ge2
$
 in an open and dense set of $\mathcal{M}\times\mathbb{R}$. By using this property, from lemmas 1.9.1, 1.9.2 and 1.9.6 in \cite{Isi95} we immediately obtain that the convergence of algorithm 5 is achieved in at most $n-1$ steps
$\blacktriangleleft$%}

%%%%%%%%%%%%%%%%%%%%%%%%%%%%%%%%%%%%%%%%
%%%%%%%%%%%%%%%%%%%%%%%%%%%%%%%%%%%%%%%%
%%%%%%%%%%%%%%%%%%%%%%%%%%%%%%%%%%%%%%%%
%%%%%%%%%%%%%%%%%%%%%%%%%%%%%%%%%%%%%%%%
%%%%%%%%%%%%%%%%%%%%%%%%%%%%%%%%%%%%%%%%

\subsection{Proof of theorem \ref{TheoremECRC}}\label{SubSectionDerivationCon}

%\proof{

\noindent We know that the first component of $\underline{x}$ (i.e., the time $t$) is not controllable and it will be not surprising to obtain this result through our analysis.

\noindent We use algorithm 3 to compute the
controllability distribution in the extended state for the system defined by (\ref{EquationChronoDynamicsOutputsTV}). We obtain:

\begin{algorithm}
\begin{algorithmic}
{\blue
  \STATE Set $\underline{\Delta}=$span$\left\{\underline{f}^1,\cdots,\underline{f}^m\right\}$ 
  \WHILE{dim$\left(\underline{\Delta}\oplus[\underline{\Delta},~\underline{f}^0]\oplus\cdots\oplus[\underline{\Delta},~\underline{f}^m]\right)>$dim$(\underline{\Delta})$}
  \STATE Set $\underline{\Delta}=\underline{\Delta}\oplus[\underline{\Delta},~\underline{f}^0]\oplus\cdots\oplus[\underline{\Delta},~\underline{f}^m]$
  \ENDWHILE}
\end{algorithmic}\label{AlgoConTVC}
\caption{Controllability distribution in the extended state}
\end{algorithm}

%\begin{algorithm}[h]
%\begin{algorithmic}
%  \State Set $k=0$
%  \State Set $\underline{\Delta}_k=span\left\{\underline{f}^1,\cdots,\underline{f}^m\right\}$ 
%  \State Set $k=k+1$
%  \State Set $\underline{\Delta}_k=\underline{\Delta}_{k-1}+[\underline{\Delta}_{k-1},~\underline{f}^0]+\sum_{i=1}^m[\underline{\Delta}_{k-1},~\underline{f}^i]$
%  \While{$dim(\underline{\Delta}_k)>dim(\underline{\Delta}_{k-1})$}
%   \State Set $k=k+1$
%  \State Set $\underline{\Delta}_k=\underline{\Delta}_{k-1}+[\underline{\Delta}_{k-1},~\underline{f}^0]+\sum_{i=1}^m[\underline{\Delta}_{k-1},~\underline{f}^i]$
%  \EndWhile
%\end{algorithmic}\label{AlgoConTVC}
%\caption{Controllability distribution in the extended state}
%\end{algorithm}

\noindent From theorem \ref{TheoremCRC} we know that algorithm 6 converges 
in an open and dense set of $\mathbb{R}\times\mathcal{M}$ and the convergent distribution is obtained {\blue in at most $(n+1)-1(=n)$ steps}.
If the convergent distribution is non singular at $(t_0,~x_0)$ and its dimension is equal to $n+1$ at $(t_0,~x_0)$, then the system is weakly locally controllable from $(t_0,~x_0)$. Conversely, 
if the system is weakly locally controllable from $(t_0,~x_0)$, the dimension of the above distribution is $n+1$ in a dense neighbourhood of $(t_0,~x_0)$.

\noindent From these results it is immediate to prove the results stated by theorem \ref{TheoremECRC}.

\noindent We have ($i=1,\cdots,m$):

\[
[\underline{f}^i,~\underline{f}^0]=
\left[
\begin{array}{cccc}
 0&0&\cdots&0\\
 \frac{\partial f^0}{\partial t}&& \frac{\partial f^0}{\partial x}&\\
\end{array}
\right]
\left[
\begin{array}{c}
 0\\
 f^i\\
\end{array}
\right]-\]
\[
\left[
\begin{array}{cccc}
 0&0&\cdots&0\\
 \frac{\partial f^i}{\partial t}&& \frac{\partial f^i}{\partial x}&\\
\end{array}
\right]
\left[
\begin{array}{c}
 1\\
 f^0\\
\end{array}
\right]=
\left[
\begin{array}{c}
 0\\
 ~[f^i,~f^0]- \frac{\partial f^i}{\partial t}\\
\end{array}
\right]
\]

\noindent and by using (\ref{EquationMixBracketTV}) we have:

\begin{equation}\label{EquationLieBracketCfif0}
[\underline{f}^i,~\underline{f}^0]=\left[
\begin{array}{c}
 0\\
 \langle f^i,~f^0\rangle\\
\end{array}
\right],~~~i=1,\cdots,m
\end{equation}

\noindent Similarly, we also obtain:

\begin{equation}\label{EquationLieBracketCfifj}
[\underline{f}^i,~\underline{f}^j]=
\left[
\begin{array}{c}
 0\\
 ~[f^i,~f^j]\\
\end{array}
\right],~~~i,j=1,\cdots,m
\end{equation}

\noindent From (\ref{EquationLieBracketCfif0}) and (\ref{EquationLieBracketCfifj}) and the initialization of the algorithm 6 ($\underline{\Delta}= \textnormal{span}\left\{\underline{f}^1,\cdots,\underline{f}^m\right\}$),  we obtain that, for every step,

\begin{equation}\label{EquationDeltaSeparation}
\underline{\Delta}=
\left[
\begin{array}{c}
 0\\
 \Delta\\
\end{array}
\right]
\end{equation}

\noindent As a result, $\Delta$ is precisely the controllability distribution we want to compute, i.e., the one that only includes the controllability properties of the system defined by (\ref{EquationObsSystemDefinitionTV}).  By using (\ref{EquationLieBracketCfif0}) and (\ref{EquationLieBracketCfifj}) we immediately obtain
 that $\Delta$ is generated by algorithm 4.

\noindent Finally,  the convergence of algorithm 6 (and consequently of algorithm 4) occurs {\blue in at most $n-1$ steps instead of $(n+1)-1(=n)$ steps}.
This is obtained by using (\ref{EquationDeltaSeparation})
and the lemmas 1.8.1, 1.8.2 and 1.8.3 in \cite{Isi95}
$\blacktriangleleft$%}

%%%%%%%%%%%%%%%%%%%%%%%%%%%%%%%%%%%%%%%%
%%%%%%%%%%%%%%%%%%%%%%%%%%%%%%%%%%%%%%%%
%%%%%%%%%%%%%%%%%%%%%%%%%%%%%%%%%%%%%%%%
%%%%%%%%%%%%%%%%%%%%%%%%%%%%%%%%%%%%%%%%
%%%%%%%%%%%%%%%%%%%%%%%%%%%%%%%%%%%%%%%%

{\blue

\section{Simple illustrative examples}\label{SectionExamples}

To illustrate the two new conditions for observability and controllability, we provide two examples. Note that they are deliberately very trivial to better figure out the main features of the two algorithms.

\subsection{Observability}

We consider the system given in (\ref{EquationObsSystemDefinitionTV}) with $m=p=1$,

\[
f\triangleq f^1= 
\left[
\begin{array}{c}
  x^1 \\
  x^2 \\
  \cdots \\
  x^n \\
\end{array}
\right], ~~~f^0= 
\left[
\begin{array}{c}
  0 \\
  0 \\
  \cdots \\
  0 \\
\end{array}
\right],~~~
h_1\triangleq h = \sum_{i=1}^n x^i t^i
\]

\noindent where, in the function $h$, $x^i$ is the $i^{th}$ component of the state and $t^i$ is $t$ to the power of $i$.
We use algorithm 2 to compute the observable codistribution.
In the following, we denote by $\Omega_k$ the codistribution returned by algorithm 2 after $k$ steps.
We obtain the following result.
$\Omega_0=$span$\{dh\}$, with:

\[
dh=tdx^1+t^2dx^2+t^3dx^3+\cdots+t^ndx^n=\sum_{i=1}^nt^idx^i
\]

\noindent We compute $\Omega_1$. We need to compute $\mathcal{L}_fh$
and $\widetilde{\mathcal{L}}_{f^0}h$. We have:

\[
\mathcal{L}_fh=h, ~~~\widetilde{\mathcal{L}}_{f^0}h=\frac{\partial h}{\partial t}=\sum_{i=1}^n ix^i t^{i-1}
\]

\noindent Hence $\Omega_1=$span$\{dh, ~d\widetilde{\mathcal{L}}_{f^0}h\}$ with:

\[
d\widetilde{\mathcal{L}}_{f^0}h=\sum_{i=1}^n i t^{i-1}dx^i
\]

\noindent Since $\dim(\Omega_1)=2>1=\dim(\Omega_0)$, we need to repeat the recursive step and compute $\Omega_2$. We need to compute the two scalar fields:

\[
\mathcal{L}_f\widetilde{\mathcal{L}}_{f^0}h=\widetilde{\mathcal{L}}_{f^0}h,
\]
\[\widetilde{\mathcal{L}}_{f^0}\widetilde{\mathcal{L}}_{f^0}h=\frac{\partial \widetilde{\mathcal{L}}_{f^0}h}{\partial t}=\sum_{i=2}^n i(i-1)x^i t^{i-2}
\]

\noindent Hence, $\Omega_2=$span$\{dh, ~d\widetilde{\mathcal{L}}_{f^0}h, ~d\widetilde{\mathcal{L}}_{f^0}\widetilde{\mathcal{L}}_{f^0}h\}$, with:

\[
d\widetilde{\mathcal{L}}_{f^0}\widetilde{\mathcal{L}}_{f^0}h=\sum_{i=1}^ni(i-1)t^{i-2}dx^i
\]

\noindent and $\dim(\Omega_2)=3>2=\dim(\Omega_1)$. By proceeding in this manner we finally obtain $\Omega_{n-1}=\Omega_{n-2}\oplus$span$\{d
\widetilde{\mathcal{L}}_{f^0}^{n-1}h\}$ and $\dim(\Omega_{n-1})=n$. We conclude that the state is weakly locally observable.

\noindent We remark that, in this driftless case, we have:

\[
\widetilde{\mathcal{L}}_{f^0}=\frac{\partial}{\partial t}
\]

\noindent Therefore, the observable codistribution is obtained by only considering the output and its time derivatives up to the $n-1$ order. In other words, the result is independent of the system input. We would obtain the weak local observability by setting $m=0$.
In addition, we would also obtain the weak local observability for the system with $m=0$ and characterized by the output

\[
h = \sum_{i=1}^n h^i(x^i) t^i
\]

\noindent where $h^1,\cdots,h^n$ are $n$ scalar functions $\mathbb{R}\rightarrow\mathbb{R}$ with nonzero derivative (the case considered above corresponds to the case $h^i(x^i)=x^i$, $\forall i$)

\noindent This result is not surprising. 
The output $h = \sum_{i=1}^n x^i t^i$ weights the components of the state in a different manner. For instance, for $n=2$, it suffices to take the output at two distinct non vanishing times, $t_1, ~t_2$, to obtain two independent equations in the two components of the state. The same holds with the output $h = \sum_{i=1}^n h^i(x^i) t^i$. In this case we first obtain $h^i(x^i),~\forall i$ and then, since the functions $h^i$ have non vanishing derivative, they can be inverted to give the components of the state.

\noindent Finally, note that the case characterized by $m=0$ and $h = \sum_{i=1}^n x^i t^i$ can be investigated by using the method in section \ref{SubSectionObsTVLinear}, for linear time-variant systems. The result that we obtain is the same.

\subsection{Controllability}

We consider the system given in (\ref{EquationObsSystemDefinitionTV}) with $m=1$,

\[
f\triangleq f^1= 
\left[
\begin{array}{c}
  x^1 \\
  x^2 \\
  \cdots \\
  x^n \\
\end{array}
\right], ~~~f^0= 
\left[
\begin{array}{c}
  t \\
  t^2 \\
  \cdots \\
  t^n \\
\end{array}
\right]
\]

\noindent We use algorithm 4 to compute the controllability distribution. As for the case of observability, we denote by $\Delta_k$ the distribution returned by algorithm 4 after $k$ steps.
We obtain the following result.

\[
\Delta_0= \textnormal{span}\{f\}= \textnormal{span}\{\left[
x^1,~x^2,~x^3,~\cdots,~x^n
\right]^T
\}
\]

\noindent We compute $\Delta_1$. We need to compute $\langle f,~f^0\rangle$. We have:

\[
\langle f,~f^0\rangle=-f^0%+[1, ~2t, ~3t^2,~\cdots,~nt^{n-1}]^T
\]

\noindent Hence:

\[
\Delta_1= \textnormal{span}\left\{f,~
f^0
\right\}
\]

\noindent Since $dim(\Delta_1)=2>1=dim(\Delta_0)$, we need to repeat the recursive step and compute $\Delta_2$. We need to compute the following vector fields: $[f, ~f]=[0,~\cdots,~0]^T$, $[f^0, ~f]=f^0$, $\langle f, ~f^0\rangle=-f^0$ and

\[
\langle f^0, f^0\rangle=-\frac{\partial f^0}{\partial t}=-[1, ~2t, ~3t^2,~\cdots,~nt^{n-1}]^T
\]

\noindent Hence:

\[
\Delta_2= \textnormal{span} \left\{f, ~f^0, \frac{\partial f^0}{\partial t}
\right\}
\] 

\noindent By proceeding in this manner we finally obtain:

\[
\Delta_{n-1}=\Delta_{n-2}\oplus  \textnormal{span}\left\{
\frac{\partial^{n-2} f^0}{\partial t^{n-2}}
\right\}
\]

\noindent and $dim(\Delta_{n-1})=n$. We conclude that the system is weakly locally controllable.

\section{Aerospace application}\label{SectionApplicationAerospace}

We consider a rocket, like a lunar module, that moves in the presence of gravity and in the absence of an atmosphere.
We assume that it is equipped with a monocular camera able to detect a point feature on the ground.
Without loss of generality, we introduce a global frame whose origin coincides with the point feature and its $z$-axis points vertically upwards. We will adopt lower-case letters to denote vectors in this frame. We define the rocket local frame as the camera frame. In addition, we assume that, with respect to this frame, the moment of inertia tensor is diagonal. In particular, by approximating the rocket with a cylinder, the vertical axis of the local frame is along the cylinder axis and, the rocket's center of gravity, belongs to this axis.
We will adopt upper-case letters to denote vectors in the local frame. 
Fig \ref{FigAppObsVISFM} illustrates our system.

%\begin{multicols}{2}
\begin{figure}[htbp]
\begin{center}
\includegraphics[width=.795\columnwidth]{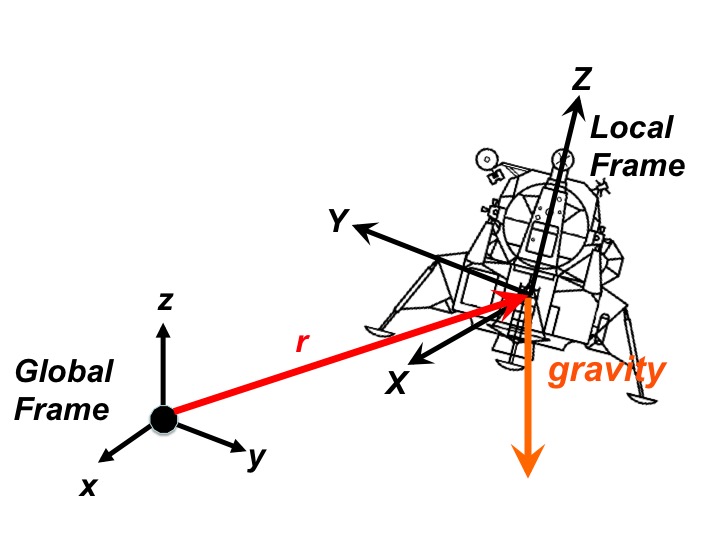}
\caption{{\blue The rocket moves in the $3D-$environment and observes a point feature at the origin by its on-board monocular camera.}} \label{FigAppObsVISFM}
\end{center}
\end{figure}

 \noindent We adopt a quaternion to represent the rocket
orientation. Indeed, even if this representation is redundant, it
is very powerful since the dynamics can be expressed in a very
easy and compact notation (see \cite{Kui99}).

%\end{multicols}

\noindent Our system is characterized by the state:

\begin{equation}\label{EquationAppObsVISFMKICalState}
x=[r_x,~r_y,~r_z,~v_x,~v_y,~v_z, ~q_t,~q_x,~q_y,~q_z,~\leo_x,~\leo_y,~\leo_z, ~g]^T
\end{equation}

\noindent where:

\begin{itemize}

\item $r=[r_x,~r_y,~r_z]^T$ is the position of the rocket in the global frame.

\item $v=[v_x,~v_y,~v_z]^T$ is the speed of the rocket in the global frame.

\item $q=q_t+q_x i+q_y j+q_z k$  is the unit quaternion that describes the rotation between the global and the local frames\footnote{A quaternion $q=q_t+q_xi+q_yj+q_zk$ is a unit quaternion if the product with its conjugate is $1$, i.e.: $qq^*=q^*q=(q_t+q_xi+q_yj+q_zk)(q_t-q_xi-q_yj-q_zk)=(q_t)^2+(q_x)^2+(q_y)^2+(q_z)^2=1$ (see \cite{Kui99}).}.

\item $\leo = \left[\leo_x~\leo_y~\leo_z \right]^T$, is the angular speed  expressed in the local frame (note that we adopted the symbol $\leo$ instead of $\Omega$ because the latter has been already used to denote the observable codistribution).

\item $g$ is the magnitude of the gravitational acceleration, which is unknown.

\end{itemize}

\noindent Additionally, we introduce the following quantities:

\begin{itemize}

%\item $A=[A_x, ~A_y, ~A_z]^T$, is the rocket acceleration in the local frame (this includes both the inertial acceleration and gravity).

\item $F$ is the magnitude of the force provided by the main engine of the rocket. The direction of this force is along the vertical axis of the local frame.

\item $\mu$ is the mass of the rocket.

\item $\mathcal{I}$ is the moment of inertia tensor of the rocket, which is diagonal $\mathcal{I}=\textnormal{diag}\{\mathcal{I}_{x},~\mathcal{I}_y,~\mathcal{I}_z\}$.

\item $\mathcal{T}\triangleq [\mathcal{T}_x, ~\mathcal{T}_y, \mathcal{T}_z]^T$ denotes the torque that acts on the rocket and which is powered by the secondary rocket engines.

\end{itemize}

\noindent

\noindent In the following, for each vector defined in the $3D$ space, the subscript $q$ will be adopted to denote the corresponding imaginary quaternion.
For instance, $\leo_q=0+\leo_x~i + \leo_y~j + \leo_z ~k$. By using the properties of the unit quaternions, we can easily obtain vectors in the global frame starting from the local frame and vice-versa. For instance, given $\leo = \left[\leo_x~\leo_y~\leo_z \right]^T$ in the local frame, we build $\leo_q=0+\leo_x~i + \leo_y~j + \leo_z ~k$, then we compute the quaternion product $\omega_q=q\leo_qq^*$. The result will be an imaginary quaternion\footnote{The product of a unit quaternion times an imaginary quaternion times the conjugate of the unit quaternion is always an imaginary quaternion.}, i.e., $\omega_q=0+\omega_x~i + \omega_y~j + \omega_z ~k$. The vector $\omega=\left[\omega_x~\omega_y~\omega_z \right]^T$ is the rocket angular speed in the global frame. Conversely, to obtain this vector in the local frame starting from $\omega$, it suffices to compute the quaternion product $\leo_q=q^*\omega_qq$.

\noindent By using this notation, the rocket acceleration generated by the main engine in the global frame is $q\frac{F}{\mu }kq^*$ and, by including the gravity, we have:

\[
\dot{v}_q =\frac{F}{\mu }~ qkq^*- g k
\]

\noindent where $k$ is the fourth fundamental quaternion unit ($k=0+0 ~i + 0~j+1~k$).

\noindent Note that, the mass $\mu$ decreases during the maneuver, due to the fuel consumption. We assume that $\mu=\mu(t)=\mu_0(1-kt)$, where $\mu_0$ is the initial mass and $k$ characterizes the consumption rate, which is constant. We denote by $T$ the time required to terminate the entire amount of fuel. Finally, we assume that $kT<<1$, meaning that the weight of the entire amount of fuel is much smaller than the weight of the rocket.
Under these assumptions, we can use the following approximation:

\begin{equation}\label{EquationApplicationMassDecrease}
\frac{F_0}{\mu }=A(t)=A_0+A_1t+A_2t^2
\end{equation}

\noindent where $A_0=F_0/\mu_0$, $	A_1=F_0k/\mu_0$ and $A_2=F_0k^2/\mu_0$.

\noindent To complete the derivation of the dynamics, we need to deal with the angular components. We start by reminding the reader the Euler's equation for the rigid body dynamics:

\begin{equation}
\mathcal{I}\dot{\leo} + \leo\wedge(\mathcal{I}\Omega)=\mathcal{T}
\end{equation}

\noindent where the symbol "$\wedge$" denotes the vector product. In our reference frame, where $\mathcal{I}$ is diagonal, we obtain the following three equations:

\begin{equation}\label{EquationAppTorque}
\left[\begin{array}{ll}
\dot{\leo}_x &= \frac{\mathcal{I}_y-\mathcal{I}_z}{\mathcal{I}_x}\leo_y\leo_z+\frac{1}{\mathcal{I}_x}\mathcal{T}_x\\

\dot{\leo}_y &= \frac{\mathcal{I}_z-\mathcal{I}_x}{\mathcal{I}_y}\leo_x\leo_z+\frac{1}{\mathcal{I}_y}\mathcal{T}_y\\

\dot{\leo}_z &= \frac{\mathcal{I}_x-\mathcal{I}_y}{\mathcal{I}_z}\leo_x\leo_y+\frac{1}{\mathcal{I}_z}\mathcal{T}_z\\
\end{array}
\right.
\end{equation}

\noindent We compute the values of the moment of inertia in our reference frame. We approximate the rocket with a cylinder. Since the center of gravity belongs to the vertical axis, we have:

\begin{equation}
\mathcal{I}_z=\frac{\mu}{2}r^2
\end{equation}

\noindent where $r$ is the radius of the approximating cylinder. The computation of the other two components, $\mathcal{I}_x$ and $\mathcal{I}_y$, is a bit more complex since, due to the fuel consumption, the center of gravity moves (we assume that it moves along the $z$-axis). We use the Parallel axis theorem \cite{Kane05}. This theorem is very simple and it allows us to compute the moment of inertia of a rigid body about any axis, given the body's moment of inertia about a parallel axis through the object's center of gravity and the perpendicular distance between the axes. Specifically:

\[
\mathcal{I}_a=\mathcal{I}^c_a+\mu d^2
\]

\noindent where $\mathcal{I}_a$ is the moment of inertia with respect to the axis $a$, $\mathcal{I}_a^c$ is the moment of inertia with respect to the axis parallel to $a$ but passing through the object's center of gravity, and $d$ is the distance between the two axes.

\noindent Because of the symmetry of the cylinder, $\mathcal{I}_x=\mathcal{I}_y\triangleq\mathcal{I}_s$. In addition, from the Parallel axis theorem we have:

\[
\mathcal{I}_s=\mathcal{I}_s^c+\mu d^2
\] 

\noindent $\mathcal{I}_s^c$ can be easily computed. For instance, for a cylinder of mass $\mu$, radius $r$ and length $h$ it is:

\[
\mathcal{I}_s^c=\frac{\mu}{12}(3r^2+h^2)
\] 

\noindent By using the same approximation given in (\ref{EquationApplicationMassDecrease}), we introduce the following approximations:

\begin{equation}\label{EquationApplicationInertiaVariation}
\begin{array}{ll}
\frac{1}{\mathcal{I}_s}&\triangleq\mathcal{S}(t)\simeq \mathcal{S}_0+\mathcal{S}_1t+\mathcal{S}_2t^2\\
\frac{1}{\mathcal{I}_z}&\triangleq\mathcal{Z}(t)\simeq \mathcal{Z}_0+\mathcal{Z}_1t+\mathcal{Z}_2t^2\\
 \frac{\mathcal{I}_s-\mathcal{I}_z}{\mathcal{I}_s}&\triangleq\mathcal{R}(t)\simeq \mathcal{R}_0+\mathcal{R}_1t+\mathcal{R}_2t^2\\
 \end{array}
\end{equation}

\noindent where all the coefficients $\mathcal{S}_i$, $\mathcal{Z}_i$ and $\mathcal{R}_i$, $i=1,2,3$ can be easily obtained from the equations above and the geometry of the rocket. By substituting in (\ref{EquationAppTorque}) we obtain:

\begin{equation}\label{EquationAppTorqueA}
\left[\begin{array}{ll}
\dot{\leo}_x &= \mathcal{R}\leo_y\leo_z+\mathcal{S}\mathcal{T}_x\\

\dot{\leo}_y &= -\mathcal{R}\leo_x\leo_z+\mathcal{S}\mathcal{T}_y\\

\dot{\leo}_z &= \mathcal{Z}\mathcal{T}_z\\
\end{array}
\right.
\end{equation}

\noindent Finally, we assume that the rocket accomplishes small movements and we can assume that the magnitude of gravity is constant. The dynamics of the state are:
\begin{equation}\label{EquationAppObsVISFMKICalDynamics}
\left[\begin{array}{ll}
\dot{r}_q &= v_q\\
\dot{v}_q &=\frac{F}{\mu } ~qkq^*- g~ k\\
\dot{q} &= \frac{1}{2}q\leo_q\\
\dot{\leo}_x &= \mathcal{R}\leo_y\leo_z+\mathcal{S}\mathcal{T}_x\\

\dot{\leo}_y &= -\mathcal{R}\leo_x\leo_z+\mathcal{S}\mathcal{T}_y\\

\dot{\leo}_z &= \mathcal{Z}\mathcal{T}_z\\
\dot{g} &=0
\end{array}
\right.
\end{equation}

\noindent We assume that the three quantities $\mathcal{T}_x$, $\mathcal{T}_y$ and $\mathcal{T}_z$ can be set by activating the secondary engines of the rocket.

\noindent The monocular camera provides the position of the point feature in the local frame, up to a scale. Let us denote this position by $P$. We have:

\begin{equation}\label{EquationAppObsVISFMKICalF}
P_q=0+P_x ~i + P_y~j+P_z~k=q^*(-r_q)q=-q^*(r_q)q
\end{equation}

\noindent Since the camera provides $P$ up to a scale, it provides the two ratios of its components: $\frac{P_x}{P_z}, ~\frac{P_y}{P_z}$. These are the outputs of our system. 
Finally, we must account for the constraint that expresses the unity of $q$. For the observability analysis,we can account for this constraint by adding the further output $(q_t)^2+(q_x)^2+(q_y)^2+(q_z)^2$. Therefore, our system is characterized by the following three outputs:

\begin{equation}\label{EquationAppObsVISFMKICalOutput}
y=h(x)=
\left[\begin{array}{c}
h_1(x)\\
h_2(x)\\
h_3(x)\\
\end{array}
\right]=
\left[\begin{array}{c}
P_x/P_z\\
P_y/P_z\\
(q_t)^2+(q_x)^2+(q_y)^2+(q_z)^2\\
\end{array}
\right]
\end{equation}

\noindent where the components $P_x$, $P_y$ and $P_z$ depend on the state $x$ through (\ref{EquationAppObsVISFMKICalF}).

\vskip.2cm 

\noindent In the following, 
we study the observability and the controllability of the state in (\ref{EquationAppObsVISFMKICalState}). Its dimension is $n=14$.
Its dynamics (given in (\ref{EquationAppObsVISFMKICalDynamics})) are characterized by $4$ inputs: $F$, $\mathcal{T}_x$, $\mathcal{T}_y$ and $\mathcal{T}_z$. On the other hand, we consider the case when the input $F$ is constant. We set its value to $F_0$. Therefore, our system is characterized by a state with dimension $n=14$, $m=3$ inputs and $p=3$  outputs (given in (\ref{EquationAppObsVISFMKICalF}) and (\ref{EquationAppObsVISFMKICalOutput})).

\noindent By comparing with equation (\ref{EquationObsSystemDefinitionTV}) we have $u_1=\mathcal{T}_x$, $u_2=\mathcal{T}_y$ and $u_3=\mathcal{T}_z$, and:

\begin{equation}\label{EquationApplicationVectors}
\small
\footnotesize
%\scriptsize
f^0(x)=\left[\begin{array}{c}
v_x\\
v_y\\
v_z\\
              2 \frac{F_0}{\mu } q_t q_y + 2 \frac{F_0}{\mu } q_x q_z\\
                 2 \frac{F_0}{\mu } q_y q_z - 2 \frac{F_0}{\mu } q_t q_x\\
 \frac{F_0}{\mu } q_t^2 - \frac{F_0}{\mu } q_x^2 - \frac{F_0}{\mu } q_y^2 + \frac{F_0}{\mu } q_z^2 - g\\
   - (\leo_x q_x)/2 - (\leo_y q_y)/2 - (\leo_z q_z)/2\\
     (\leo_x q_t)/2 - (\leo_y q_z)/2 + (\leo_z q_y)/2\\
     (\leo_y q_t)/2 + (\leo_x q_z)/2 - (\leo_z q_x)/2\\
     (\leo_z q_t)/2 - (\leo_x q_y)/2 + (\leo_y q_x)/2\\
\mathcal{R}\leo_y\leo_z\\
-\mathcal{R}\leo_x\leo_z\\
0\\
0\\
\end{array}
\right],~
f^1(x)=\left[\begin{array}{c}
0\\
0\\
0\\
0\\
0\\
0\\
0\\
0\\
0\\
0\\
\mathcal{S}\\
0\\
0\\
0\\
\end{array}
\right],~f^2(x)=\left[\begin{array}{c}
0\\
0\\
0\\
0\\
0\\
0\\
0\\
0\\
0\\
0\\
0\\
\mathcal{S}\\
0\\
0\\
\end{array}
\right],~f^3(x)=\left[\begin{array}{c}
0\\
0\\
0\\
0\\
0\\
0\\
0\\
0\\
0\\
0\\
0\\
0\\
\mathcal{Z}\\
0\\
\end{array}
\right]
\end{equation}

\noindent The dynamics explicitly depend over time. In particular, 
$\frac{F_0}{\mu }$ is given in (\ref{EquationApplicationMassDecrease}) and $\mathcal{R}, ~\mathcal{S}, ~\mathcal{Z}$ are given in (\ref{EquationApplicationInertiaVariation}).

\subsection{Observability}
In order to obtain the observable codistribution, we need to run algorithm 2.
As in section \ref{SectionExamples}, we denote by $\Omega_k$ the codistribution returned by algorithm 2 after $k$ steps.
We obtain, at the initialization:

\[
\Omega_0=\textnormal{span}\{dh_1,~dh_2,~dh_3\}
\]

\noindent Its dimension is $3$. We compute $\Omega_1$. We obtain:

\[
\Omega_1=\Omega_0\oplus \textnormal{span}\{\widetilde{\mathcal{L}}_{f^0} dh_1,~\widetilde{\mathcal{L}}_{f^0} dh_2\}
\]

\noindent Note that the functions $h_1$ and $h_2$ do not explicitly depend on time and $\widetilde{\mathcal{L}}_{f^0}dh_i=\mathcal{L}_{f^0}dh_i$, $i=1,~2$.
We obtain $\dim(\Omega_1)=5>3=\dim(\Omega_0)$. Hence we need to compute $\Omega_2$. We obtain:

\[
\Omega_2=\Omega_1\oplus \textnormal{span}\{\widetilde{\mathcal{L}}^2_{f^0} dh_1,~\widetilde{\mathcal{L}}^2_{f^0} dh_2\}
\]

\noindent Again, the functions $\widetilde{\mathcal{L}}_{f^0}h_1$ and $\widetilde{\mathcal{L}}_{f^0}h_2$ do not explicitly depend on time and $\widetilde{\mathcal{L}}^2_{f^0}dh_i=\mathcal{L}^2_{f^0}dh_i$, $i=1,~2$.
We obtain  $\dim(\Omega_2)=7>5=\dim(\Omega_1)$. Hence we need to compute 
$\Omega_3$. We obtain:

\[
\Omega_3=\Omega_2\oplus \textnormal{span}\{\widetilde{\mathcal{L}}^3_{f^0} dh_1,~\widetilde{\mathcal{L}}^3_{f^0} dh_2, ~\mathcal{L}_{f^1}\mathcal{L}^2_{f^0} dh_1, ~\mathcal{L}_{f^2}\mathcal{L}^2_{f^0} dh_1, ~\mathcal{L}_{f^1}\mathcal{L}^2_{f^0} dh_2\}
\]

\noindent This time, the two functions $\widetilde{\mathcal{L}}^2_{f^0}h_1$ and $\widetilde{\mathcal{L}}^2_{f^0}h_2$ explicitly depend on time and

\[
\widetilde{\mathcal{L}}^3_{f^0}dh_i=\widetilde{\mathcal{L}}_{f^0}\mathcal{L}^2_{f^0}dh_i=
\mathcal{L}^3_{f^0}dh_i+\frac{\partial}{\partial t}\mathcal{L}^2_{f^0}dh_i
\neq\mathcal{L}^3_{f^0}dh_i
\]

\noindent $i=1,~2$.
We obtain  $\dim(\Omega_3)=12>7=\dim(\Omega_2)$. Hence, we need to compute 
$\Omega_4$. We obtain:

\[
\Omega_4=\Omega_3\oplus
 \textnormal{span}\{\mathcal{L}_{f^2}\widetilde{\mathcal{L}}^3_{f^0} dh_1\}
\]

\noindent and $\dim(\Omega_4)=13>12=\dim(\Omega_3)$. Hence, we need to compute 
$\Omega_5$. We obtain:

\[
\Omega_5=\Omega_4
\]

\noindent This means that algorithm 2 converges at its fourth step and the observable codistribution is $\Omega=\Omega_4$.
As a result, its dimension is $13<n=14$ and the state is not observable. By computing the orthogonal distribution we obtain: 

\[
\Omega^{\bot}=
 \textnormal{span}\left\{
\left[\begin{array}{c}
-r_y\\
r_x\\
0\\
-v_y\\
v_x\\
0\\
-q_z/2\\
-q_y/2\\
q_x/2\\
q_t/2\\
0\\
0\\
0\\
0\\
\end{array}
\right]
\right\}
\]

\noindent The infinitesimal transformation discussed in \cite{TRO10} becomes, for the specific case:

\[
\left[
\begin{array}{c}
    r_x \\
    r_y\\
    r_z\\
\end{array}
\right]\rightarrow \left[
\begin{array}{c}
    r_x \\
    r_y\\
    r_z\\
\end{array}
\right] + \epsilon \left[
\begin{array}{c}
    -r_y \\
    r_x\\
    0\\
\end{array}
\right]\simeq
\left[\begin{array}{ccc}
 \cos\epsilon & -\sin\epsilon & 0 \\
\sin\epsilon & \cos\epsilon & 0\\
0& 0& 1\\
\end{array}\right]
\left[
\begin{array}{c}
    r_x \\
    r_y\\
    r_z\\
\end{array}
\right]
\]

\[
\left[
\begin{array}{c}
    v_x \\
    v_y\\
    v_z\\
\end{array}
\right]\rightarrow \left[
\begin{array}{c}
    v_x \\
    v_y\\
    v_z\\
\end{array}
\right] + \epsilon \left[
\begin{array}{c}
    -v_y \\
    v_x\\
    0\\
\end{array}
\right]\simeq
\left[\begin{array}{ccc}
 \cos\epsilon & -\sin\epsilon & 0 \\
\sin\epsilon & \cos\epsilon & 0\\
0& 0& 1\\
\end{array}\right]
\left[
\begin{array}{c}
    v_x \\
    v_y\\
    v_z\\
\end{array}
\right]
\]

\[
\left[
\begin{array}{c}
    q_t \\
    q_x \\
    q_y\\
    q_z\\
\end{array}
\right]\rightarrow \left[
\begin{array}{c}
    q_t \\
    q_x \\
    q_y\\
    q_z\\
\end{array}
\right] + \frac{\epsilon}{2} \left[
\begin{array}{c}
    -q_z \\
    -q_y \\
    q_x\\
    q_t\\
\end{array}
\right]
\]

\noindent that is an infinitesimal rotation about the vertical axis (regarding the quaternion, this can be verified starting from the last equation in (\ref{EquationAppObsVISFMKICalDynamics}) that provides $\dot{q}=\frac{1}{2}\omega q$ (where $\omega$ is the angular speed in the global frame) and by setting $\omega_x=\omega_y=0$ and $\omega_z dt=\epsilon$).

\vskip.2cm

\noindent We want to describe our system with an observable state. We have many choices. Since we want to achieve a local observable subsystem, we are interested in introducing very simple observable functions that generate the observable codistribution. First of all, we remark that the distance of the point feature is an observable function. This can be verified by checking that the gradient of the function $r_x^2+r_y^2+r_z^2$ belongs to $\Omega$ (actually, we even do not need to check this: it suffices to remark that the scale is rotation invariant and, consequently, it satisfies the above symmetry and it is observable). Now, if the scale is observable, by combining this knowledge with the two outputs $h_1=\frac{P_x}{P_z}$ and $h_2=\frac{P_y}{P_z}$ (and by knowing that $P_x^2+P_y^2+P_z^2=r_x^2+r_y^2+r_z^2$), we can select the following three observable functions:

\[
P_x, ~~P_y, ~~P_z
\]

\noindent In a similar manner, we can also select the following three observable functions:

\[
V_x, ~~V_y, ~~V_z
\]

\noindent which are the components of the body speed expressed in the local frame. Finally, regarding the orientation, we know that the yaw, which characterizes a rotation about the vertical axis, is unobservable. Conversely, the roll and the pitch angles are observable. We denote the roll and pitch with $\psi_R$ and $\psi_P$, respectively.
Hence, we introduce the following observable state:

\begin{equation}\label{EquationAppObsVISFMKICalStateObservable}
X=[P_x,~P_y,~P_z,~V_x,~V_y,~V_z, ~\psi_R, ~\psi_P, ~\leo_x, ~\leo_y, ~\leo_z, ~g]^T
\end{equation}

\noindent These components can be expressed in terms of the components of the original state in (\ref{EquationAppObsVISFMKICalState}) as follows. The first three components are given by (\ref{EquationAppObsVISFMKICalF}). The following five are:

\[
V_q=0+V_x ~i + V_y~j+V_z~k=q^*v_qq,~
\psi_R=\arctan\left(
2\frac{q_tq_x+q_yq_z}{1-2(q_x^2+q_y^2)}
\right),~
\psi_P=\arcsin\left(
2(q_tq_y-q_xq_z)
\right)
\]

\noindent From the above equations, we obtain the following description of the local observable subsystem:

\begin{equation}\label{EquationAppObsVISFMLocDec}
\left[\begin{array}{ll}
\dot{P} &=-\leo \wedge P -V\\
\dot{V} &=-\leo \wedge V +A - G\\
\dot{\psi}_R &= \leo_x + \leo_y \tan\psi_P\sin\psi_R + \leo_z \tan\psi_P\cos\psi_R\\
\dot{\psi}_P &=\leo_y \cos\psi_R - \leo_z \sin\psi_R\\
\dot{\leo}_x &= \mathcal{R}\leo_y\leo_z+\mathcal{S}\mathcal{T}_x\\

\dot{\leo}_y &= -\mathcal{R}\leo_x\leo_z+\mathcal{S}\mathcal{T}_y\\

\dot{\leo}_z &= \mathcal{Z}\mathcal{T}_z\\
\dot{g} &=0\\
y &= [P_x/P_z, ~P_y/P_z]^T\\
\end{array}
\right.
\end{equation}

\noindent where $A\triangleq[0,~0,~\frac{F_0}{\mu }]^T$ and $G$ is the gravity in the local frame and it only depends on the roll and pitch angles:

\begin{equation}\label{EquationAppObsVISFMG}
G=
g\left[\begin{array}{c}
\sin\psi_P\\
-\cos\psi_P\sin\psi_R\\
-\cos\psi_P\cos\psi_R\\
\end{array}
\right]
\end{equation}

\noindent We conclude this section with the following remarks:

\begin{itemize}

\item The same observable state, and the same local subsystem given by (\ref{EquationAppObsVISFMLocDec}), still hold by considering a varying $F$. Indeed, removing the constraint $F=F_0$ cannot decrease the observability properties. On the other hand, the system remains invariant under rotations about the gravity and the absolute yaw remains unobservable.

\item In the case of multiple features, the observability properties remain the same. The yaw angle remains unobservable. To prove this it is unnecessary to repeat the computation of the observable codistribution. Since the gravity is invariant to the yaw, the system maintains the same continuous symmetry that describes a rotation about the vertical axis.
In presence of $M$ point features, a local observable subsystem is given by (\ref{EquationAppObsVISFMLocDec}), where the first equation $\dot{P} =-\leo \wedge P -V$ must be replaced by the  $M$ equations $\dot{P}^j=-\leo \wedge P^j -V$, $j=1,\cdots,M$ and the last equation $y = [P_x/P_z, ~P_y/P_z]^T$, by $y = [P_x^1/P_z^1, ~P_y^1/P_z^1, ~\cdots, P_x^M/P_z^M, ~P_y^M/P_z^M]^T$.

\item We obtain for our system a single symmetry which expresses the unobservability of the absolute yaw angle. On the other hand, if we do not a priori know the position of the point feature in the global frame, we have three further symmetries which express the unobservability of the absolute position of the point feature (i.e., its position in this global frame). In our analysis, by introducing a global frame whose origin coincides with the point feature, we are implicitly assuming that the position of the point feature is a priori known.

\item The equations given in (\ref{EquationAppObsVISFMLocDec}) could be used for a practical implementation, (e.g., to implement an extended Kalman filter). In this case, we recommend to use a different output. The camera provides $P$ up to a scale. Instead of the two ratios, $\frac{P_x}{P_z}$ and $\frac{P_y}{P_z}$, which are singular when $P_z=0$, it is better to introduce two angles (e.g., by setting $[P_x,~P_y,~P_z]^T=|P|[\cos\alpha_1\cos\alpha_2, ~\cos\alpha_1\sin\alpha_2, ~\sin\alpha_1]^T$).

\end{itemize}

\subsection{Controllability}
 
%We study the controllability of the state in (\ref{EquationAppObsVISFMKICalState}). As for the observability, we consider the case when the input $F$ is constant ($F=F_0$). Therefore, our system is characterized by a state with dimension $n=14$, and $m=3$ inputs. By comparing with equation (\ref{EquationObsSystemDefinitionTV}) we have $u_1=\mathcal{T}_x$, $u_2=\mathcal{T}_y$ and $u_3=\mathcal{T}_z$, and the four vector fields $f^0(x)$, $f^1(x)$, $f^2(x)$ and $f^3(x)$, are given in (\ref{EquationApplicationVectors}).
%
%
%\noindent The dynamics explicitly depend over time. In particular, 
%$\frac{F_0}{\mu }$ is given in (\ref{EquationApplicationMassDecrease}) and $\mathcal{R}, ~\mathcal{S}, ~\mathcal{Z}$ are given in (\ref{EquationApplicationInertiaVariation}).
%
%

 In order to obtain the controllable distribution, we need to run algorithm 4.
As in section \ref{SectionExamples}, we denote by $\Delta_k$ the distribution returned by algorithm 4 after $k$ steps.
We obtain, at the initialization:

\[
\Delta_0=\textnormal{span}\{f^1,~f^2,~f^3\}
\]

\noindent Its dimension is $3$. We compute $\Delta_1$. We obtain:

\[
\Delta_1=\Delta_0\oplus \textnormal{span}\{
\langle
f^1,~f^0
\rangle,
\langle
f^2,~f^0
\rangle,
\langle
f^3,~f^0
\rangle
\}
\]

\noindent and $\dim(\Delta_1)=6>3=\dim(\Delta_0)$. Hence we need to compute $\Delta_2$. We obtain:

\[
\Delta_2=\Delta_1\oplus \textnormal{span}\{
\langle
\langle
f^1,~f^0
\rangle,~f^0
\rangle,
\langle
\langle
f^2,~f^0
\rangle,~f^0
\rangle
\}
\]

\noindent and $\dim(\Delta_2)=8>6=\dim(\Delta_1)$. Hence we need to compute $\Delta_3$. We obtain:

\[
\Delta_3=\Delta_2\oplus \textnormal{span}\{
\langle
\langle
\langle
f^1,~f^0
\rangle,~f^0
\rangle,~f^0
\rangle,
\langle
\langle
\langle
f^2,~f^0
\rangle,~f^0
\rangle,~f^0
\rangle
\}
\]

\noindent and $\dim(\Delta_3)=10>8=\dim(\Delta_2)$. Hence we need to compute $\Delta_4$. We obtain:

\[
\Delta_4=\Delta_3\oplus \textnormal{span}\{
[
\langle
\langle
\langle
f^1,~f^0
\rangle,~f^0
\rangle,~f^0
\rangle,~f^1
],
\langle
\langle
\langle
\langle
f^1,~f^0
\rangle,~f^0
\rangle,~f^0
\rangle,~f^0
\rangle
\}
\]

\noindent and $\dim(\Delta_4)=12>9=\dim(\Delta_3)$. Hence we need to compute $\Delta_5$. We obtain:

\[
\Delta_5=\Delta_4
\]

\noindent This means that algorithm 4 converges at its fourth step and the controllable distribution is $\Delta=\Delta_4$.
By computing the orthogonal distribution we obtain: 

\[
\Delta^{\bot}=
 \textnormal{span}\left\{dg, ~~q_tdq_t+q_xdq_x+q_ydq_y +q_zdq_z
\right\}
\]

\noindent The physical meaning of these two covectors is clear. The former, $dg$, means that we cannot control the last state component (i.e., $g$). This is obvious since we cannot modify the magnitude of the gravity. The latter, $q_tdq_t+q_xdq_x+q_ydq_y +q_zdq_z$, simply expresses the fact that we cannot modify the norm of the quaternion. This is obvious since the quaternion $q$ must characterize a rotation and must be a unit quaternion. Hence, we conclude that our system is weakly locally controllable.

\noindent Note that, even if we can control our system, the unobservability of the yaw angle results in the impossibility of setting the yaw to a desired value. However, in most of cases, this is unnecessary.

\section{Observability during take off and landing}\label{SectionTakeOff}

We investigate the observability properties of our system when the lunar module undertakes two very important maneuvers: take off and landing. During both these maneuvers, the lunar module attitude remains constant. Under these conditions, to obtain the observability properties, we can 
characterize our system by the following reduced state:

\begin{equation}\label{EquationAppObsVISFMKICalStateTO}
x=[r_x,~r_y,~r_z,~v_x,~v_y,~v_z, ~q_t,~q_x,~q_y,~q_z, ~g]^T
\end{equation}

\noindent Its dynamics are:
\begin{equation}\label{EquationAppObsVISFMKICalDynamicsTO}
\left[\begin{array}{ll}
\dot{r}_q &= v_q\\
\dot{v}_q &=\frac{F}{\mu } ~qkq^*- g~ k\\
\dot{q} &= 0\\
\dot{g} &=0
\end{array}
\right.
\end{equation}

\noindent Note that, this characterization is only used to obtain the observability properties. The attitude of the lunar module is actually maintained constant by activating the secondary engines. These engines provide suitable values of $\mathcal{T}_x$, $\mathcal{T}_y$ and $\mathcal{T}_z$ in the dynamics given in (\ref{EquationAppObsVISFMKICalDynamics}), in order to maintain $\dot{q}=0$, as in (\ref{EquationAppObsVISFMKICalDynamicsTO}). In addition, we should set $q=1$ and remove it from the state. However, in order to maintain $q=1$ we need to check that $q$ is observable. We actually know that this is not the case, since, even without restrictions on the rocket motion, the absolute yaw is unobservable. Hence, what we can hope, is to maintain $q(t)=\cos\frac{\psi_Y(t)}{2}+\sin\frac{\psi_Y(t)}{2}~k$, where $\psi_Y$ is the yaw angle.

\noindent In the following, we investigate two distinct scenarios. In both, we always set $F=F_0$.
The former is characterized by also a constant value of the mass ($\mu=\mu_0$). In this case, $\frac{F_0}{\mu }$ is independent of time and we set $\frac{F_0}{\mu_0}=A_0$. In the latter, the mass decreases during the operation, due to the fuel consumption. 
This is precisely the same case analyzed in section \ref{SectionApplicationAerospace} and
$\frac{F_0}{\mu }$ is given by equation
(\ref{EquationApplicationMassDecrease}).

\subsection{First scenario: observability with constant mass}

\noindent By comparing equation (\ref{EquationAppObsVISFMKICalDynamicsTO}) with (\ref{EquationObsSystemDefinitionTV}) we have no input and:

%\begin{equation}\label{EquationApplicationVectorsTO}
%%\small
%%\ootnotesize
%%\scriptsize
%f^0(x)=\left[\begin{array}{c}
%v_x\\
%v_y\\
%v_z\\
%              2 A_0 q_t q_y + 2A_0 q_x q_z\\
%                 2 A_0 q_y q_z - 2A_0 q_t q_x\\
%A_0 q_t^2 - A_0 q_x^2 - A_0 q_y^2 + A_0 q_z^2 - g\\
%   - (\leo_x q_x)/2 - (\leo_y q_y)/2 - (\leo_z q_z)/2\\
%     (\leo_x q_t)/2 - (\leo_y q_z)/2 + (\leo_z q_y)/2\\
%     (\leo_y q_t)/2 + (\leo_x q_z)/2 - (\leo_z q_x)/2\\
%     (\leo_z q_t)/2 - (\leo_x q_y)/2 + (\leo_y q_x)/2\\
%0\\
%\end{array}
%\right]
%\end{equation}
%

\begin{equation}\label{EquationApplicationVectorsTO}
%\small
%\ootnotesize
%\scriptsize
f^0(x)=\left[\begin{array}{c}
v_x\\
v_y\\
v_z\\
              2 A_0 q_t q_y + 2A_0 q_x q_z\\
                 2 A_0 q_y q_z - 2A_0 q_t q_x\\
A_0 q_t^2 - A_0 q_x^2 - A_0 q_y^2 + A_0 q_z^2 - g\\
0\\
0\\
0\\
0\\
0\\
\end{array}
\right]
\end{equation}

\noindent In addition, the system has $p=3$  outputs (given in (\ref{EquationAppObsVISFMKICalF}) and (\ref{EquationAppObsVISFMKICalOutput})).

\noindent The dynamics do not explicitly depend over time. 
In order to obtain the observable codistribution, we need to run algorithm 1.
As in section \ref{SectionExamples}, we denote by $\Omega_k$ the codistribution returned by algorithm 1 after $k$ steps.
We obtain, at the initialization:

\[
\Omega_0=\textnormal{span}\{dh_1,~dh_2,~dh_3\}
\]

\noindent Its dimension is $3$. We compute $\Omega_1$. We obtain:

\[
\Omega_1=\Omega_0\oplus \textnormal{span}\{\mathcal{L}_{f^0} dh_1,~\mathcal{L}_{f^0}dh_2\}
\]

\noindent We obtain $\dim(\Omega_1)=5>3=\dim(\Omega_0)$. Hence we need to compute $\Omega_2$. We obtain:

\[
\Omega_2=\Omega_1\oplus \textnormal{span}\{\mathcal{L}^2_{f^0} dh_1,~\mathcal{L}^2_{f^0} dh_2\}
\]

\noindent We obtain  $\dim(\Omega_2)=7>5=\dim(\Omega_1)$. Hence we need to compute 
$\Omega_3$. We obtain:

\[
\Omega_3=\Omega_2\oplus \textnormal{span}\{\mathcal{L}^3_{f^0} dh_1,~\mathcal{L}^3_{f^0} dh_2\}
\]

\noindent We obtain  $\dim(\Omega_3)=9>7=\dim(\Omega_2)$. Hence, we need to compute 
$\Omega_4$. We obtain:

%\[
%\Omega_4=\Omega_3\oplus
% \textnormal{span}\{\widetilde{\mathcal{L}}^4_{f^0} dh_1\}
%\]
%
\[
\Omega_4=\Omega_3
\]

\noindent This means that algorithm 1 converges at its third step and the observable codistribution is $\Omega=\Omega_3$.
As a result, its dimension is $9<n=11$ and the state is not observable. By computing the orthogonal distribution we obtain:

\[
\Omega^{\bot}=
 \textnormal{span}\left\{
 \left[\begin{array}{c}
-r_y\\
r_x\\
0\\
-v_y\\
v_x\\
0\\
-q_z/2\\
-q_y/2\\
q_x/2\\
q_t/2\\
0\\
\end{array}
\right],~~~
\left[
\begin{array}{c}
   A_0r_z Q_x+gr_x\\
   A_0r_z  Q_y+gr_y\\
 A_0(-   r_x Q_x -    r_y Q_y)+gr_z\\
  A_0 v_z Q_x+gv_x\\
   A_0v_z  Q_y+gv_y\\
 A_0(-   v_x Q_x -    v_y Q_y)+gv_z\\
                     A_0(-  q_t q_x^2 -   q_t q_y^2 )\\
                   A_0(q_t^2 q_x +    q_x q_z^2 )\\
              A_0 (q_t^2 q_y +    q_y q_z^2 )\\
 A_0(-q_z q_x^2 - q_zq_y^2)\\
  A_0g (-   q_t^2 +   q_x^2 +   q_y^2 -   q_z^2)+g^2\\
\end{array}
\right]
\right\}
\]

\noindent where $Q_x=2(q_t q_y +    q_z q_x)$, $Q_y=2(-q_t q_x +    q_z q_y)$. The former characterizes the invariance of the system under rotations about the gravity axis. This expresses the already known unobservability of the yaw. 
\noindent The physical meaning of the latter generator is more complex. 
It expresses the unobservability of the absolute scale. To visualize this, it is better to consider the case when $q=1$ (i.e., $q_t=1$ and $q_x=q_y=q_z=0$). The second generator becomes (up to a factor):

\[
\left[
\begin{array}{c}
   r_x\\
  r_y\\
 r_z\\
 v_x\\
 v_y\\
v_z\\
0\\
0\\
0\\
0\\
g-A_0\\
\end{array}
\right]
\]

\noindent Regarding the first 6 components of the state (i.e., the rocket position and speed in the global frame), it characterizes precisely a scale transformation. Restricted to the first six state components, the infinitesimal transformation discussed in \cite{TRO10} becomes, for the specific case:

\[
\left[
\begin{array}{c}
    r_x \\
    r_y\\
    r_z\\
    v_x \\
    v_y\\
    v_z\\
\end{array}
\right]\rightarrow \left[
\begin{array}{c}
    r_x \\
    r_y\\
    r_z\\
    v_x \\
    v_y\\
    v_z\\
\end{array}
\right] + \epsilon \left[
\begin{array}{c}
    r_x \\
    r_y\\
    r_z\\
    v_x \\
    v_y\\
    v_z\\
\end{array}
\right]=(1+\epsilon)
\left[
\begin{array}{c}
    r_x \\
    r_y\\
    r_z\\
    v_x \\
    v_y\\
    v_z\\
\end{array}
\right],
\]

\noindent which is an infinitesimal scale transform, for the rocket position and speed.
This means that we need to equip our rocket with a further sensor able to provide the scale (e.g., a laser range finder).

\subsection{Second scenario: observability with varying mass}

In this case the dynamics explicitly depend over time. In particular, 
the drift ($f^0$) is given in (\ref{EquationApplicationVectorsTO}) where, instead of $A_0$, we must substitute the expression of $\frac{F_0}{\mu }$ given in (\ref{EquationApplicationMassDecrease}). In order to obtain the observable codistribution, we need to run algorithm 2.
As in section \ref{SectionExamples}, we denote by $\Omega_k$ the codistribution returned by algorithm 2 after $k$ steps.
We obtain, at the initialization:

\[
\Omega_0=\textnormal{span}\{dh_1,~dh_2,~dh_3\}
\]

\noindent Its dimension is $3$. We compute $\Omega_1$. We obtain:

\[
\Omega_1=\Omega_0\oplus \textnormal{span}\{\widetilde{\mathcal{L}}_{f^0} dh_1,~\widetilde{\mathcal{L}}_{f^0} dh_2\}
\]

\noindent We obtain $\dim(\Omega_1)=5>3=\dim(\Omega_0)$. Hence we need to compute $\Omega_2$. We obtain:

\[
\Omega_2=\Omega_1\oplus \textnormal{span}\{\widetilde{\mathcal{L}}^2_{f^0} dh_1,~\widetilde{\mathcal{L}}^2_{f^0} dh_2\}
\]

\noindent We obtain  $\dim(\Omega_2)=7>5=\dim(\Omega_1)$. Hence, we need to compute 
$\Omega_3$. We obtain:

\[
\Omega_3=\Omega_2\oplus \textnormal{span}\{\widetilde{\mathcal{L}}^3_{f^0} dh_1,~\widetilde{\mathcal{L}}^3_{f^0} dh_2\}
\]

\noindent We obtain  $\dim(\Omega_3)=9>7=\dim(\Omega_2)$. Hence, we need to compute 
$\Omega_4$. We obtain:

\[
\Omega_4=\Omega_3\oplus
 \textnormal{span}\{\widetilde{\mathcal{L}}^4_{f^0} dh_1\}
\]

\noindent and $\dim(\Omega_4)=10>9=\dim(\Omega_3)$. Hence, we need to compute 
$\Omega_5$. We obtain:

\[
\Omega_5=\Omega_4
\]

\noindent This means that algorithm 2 converges at its fourth step and the observable codistribution is $\Omega=\Omega_4$. We have: 

\[
\Omega^{\bot}=
 \textnormal{span}\left\{
\left[\begin{array}{c}
-r_y\\
r_x\\
0\\
-v_y\\
v_x\\
0\\
-q_z/2\\
-q_y/2\\
q_x/2\\
q_t/2\\
0\\
\end{array}
\right]
\right\}
\]

\noindent namely, with respect to the case analyzed in the previous section, we gain the information on the absolute scale but we still have the system invariance with respect to a rotation about the gravity axis, as expected.

\vskip.2cm

%\begin{Rm}
%The fact that the absolute scale becomes observable in the second scenario (varying mass), can be explained by using the results obtained in \cite{TRO12} and \cite{IJCV14}. These works prove that a motion characterized by a constant inertial acceleration is singular. Specifically, in \cite{TRO12}, we proved that, when the magnitude of the gravity is unknown (as in our case) and the motion is characterized by a constant acceleration, the absolute scale is unobservable\footnote{In \cite{IJCV14} we also proved that, if the magnitude of the gravity is known, the scale becomes observable but the problem has two distinct solutions instead of a unique solution. To have a unique solution, the inertial acceleration must be time varying.}. Since we are assuming that the force provided by the main engine is constant ($F=F_0$) and the attitude is constant, when the mass does not variate (first scenario), the inertial acceleration that characterizes the motion is also constant. On the other hand, if the mass decreases during the maneuver, the inertial acceleration is not constant. 
%\end{Rm}

\begin{Rm}
The fact that the absolute scale becomes observable in the second scenario (varying mass), can be explained by using the results obtained in \cite{TRO12}. This work proves that, when the magnitude of the gravity is unknown (as in our case) and the motion is characterized by a constant acceleration, the absolute scale is unobservable. Since we are assuming that the force provided by the main rocket engine is constant ($F=F_0$) and the attitude is constant, when the mass does not variate (first scenario), the inertial acceleration that characterizes the motion is also constant. On the other hand, if the mass decreases during the maneuver, the inertial acceleration is not constant. 
\end{Rm}

%\[
%\Omega^{\bot}=
% \textnormal{span}\left\{
%2A_0
%\left[
%\begin{array}{c}
%   r_z q_t q_y +    q_z r_z  q_x\\
% -    r_z q_t q_x +    q_z r_z q_y \\
% -   q_t q_y r_x -    q_x q_z  r_x +    q_t q_x r_y -    q_y q_z r_y\\
%      v_z q_t q_y +    q_z v_z q_x \\
%   -    v_z q_t q_x +    q_z v_z q_y \\
% -   q_t q_y v_x +    q_t q_x v_y -    q_x q_z v_x -    q_y q_z v_y\\
%                     (-  q_t q_x^2 -   q_t q_y^2 )/2\\
%                   (q_t^2 q_x +    q_x q_z^2 )/2\\
%               (q_t^2 q_y +    q_y q_z^2 )/2\\
% (-q_z q_x^2 - q_zq_y^2)/2\\
%  g (-   q_t^2 +   q_x^2 +   q_y^2 -   q_z^2)/2\\
%
%\end{array}
%\right]
%+g
%\left[
%\begin{array}{c}
% r_x \\
% r_y \\
% r_z \\
%v_x \\
%v_y\\
%v_z\\
% 0\\
% 0\\
% 0\\
% 0\\
%g\\
%\end{array}
%\right]
%\right\}
%\]

}
\section{Conclusion}\label{SectionConclusion}

\noindent This paper extended the analytic conditions to check the weak local observability and the weak local controllability to the case of time-varying nonlinear systems. In other words, it extended the observability rank condition and the controllability rank condition to the time-varying case.

%Previous conditions only accounted for nonlinear systems that do not explicitly depend on time, or, for systems characterized by an explicit time dependence, they only accounted for the linear case.

\noindent The paper showed that these two conditions coincide with the well known conditions in the two simpler cases of time-varying linear systems and  time-invariant nonlinear systems. 

{\blue
\noindent The two new conditions were illustrated by discussing simple examples and by also studying the observability and the controllability properties of a  a lunar module. In particular, we analyzed this system under the constraint that the main rocket engine delivers constant power.
For this system, the dynamics exhibit an explicit time-dependence due to its the variation of the weight and the variation of the moment of inertia during the maneuvers. These variations are a consequence of the fuel consumption. To study the observability and the controllability properties of this system, the extended observability rank condition and the extended controllability rank condition introduced by this paper were requested.
We obtained that the state is weakly locally observable.
Regarding the observability, we obtained that the state is weakly locally observable with the exception of the yaw angle. In addition, we showed that, during the take off and landing, the observability of the absolute scale is a consequence of the mass variation.

\appendix

\section{Computation of the observable codistribution}\label{AppendixCod}

We provide a property that plays a key role for the implementation of algorithm 2.
We denote by $\Omega_k$ the codistribution returned by algorithm 2 at the $k^{th}$ step.
This property states that, at each step $k$, it suffices to compute the generators of $\Omega_k$ by performing simple operations (Lie derivative and time derivative) on the generators of $\Omega_{k-1}$. Additionally, these generators are differentials of scalar fields. 

\noindent We remind the reader the following property (that is the analogous that holds for algorithm 1):

\begin{pr}\label{PropertyObsOld}
Let us consider a nonsingular codistribution $\Omega$ spanned by the covectors $\omega_1,\cdots,\omega_s$ and a smooth vector field $f$. 
We have:
\[
\Omega\oplus\mathcal{L}_f\Omega= \textnormal{span}\{\omega_1,\cdots,\omega_s,\mathcal{L}_f\omega_1,\cdots,\mathcal{L}_f\omega_s\}
\]
\end{pr}

\proof{The reader is addressed to \cite{Isi95}, Remark 1.6.7
$\blacktriangleleft$}

\noindent We need to extend the above property to the case when, instead of the operator $\mathcal{L}_f$, we consider the operator defined in (\ref{EquationMixOperatorTV}).

\noindent We have the following new property.

\begin{pr}\label{PropertyObsNew}
Let us consider a nonsingular codistribution $\Omega$ spanned by the covectors $\omega_1,\cdots,\omega_s$ and a smooth vector field $f$. 
We have:
\[
\Omega\oplus\widetilde{\mathcal{L}}_f\Omega= \textnormal{span}\{\omega_1,\cdots,\omega_s,\widetilde{\mathcal{L}}_f\omega_1,\cdots,\widetilde{\mathcal{L}}_f\omega_s\}
\]
\end{pr}

\proof{Obviously, $\textnormal{span}\{\omega_1,\cdots,\omega_s,\widetilde{\mathcal{L}}_f\omega_1,\cdots,\widetilde{\mathcal{L}}_f\omega_s\}$
is included in 
$\Omega\oplus\widetilde{\mathcal{L}}_f\Omega$.

\noindent Let us prove the vice-versa. Let us consider a generic covector $\lambda\in\Omega$. We have:

\[
\lambda=\sum_{i=1}^sc_i\omega_i
\]

\noindent We have:

\[
\widetilde{\mathcal{L}}_f\lambda=
\mathcal{L}_f \left(\sum_{i=1}^sc_i\omega_i\right)+
\frac{\partial}{\partial t} \left(\sum_{i=1}^sc_i\omega_i\right)=
\sum_{i=1}^s  \left(\mathcal{L}_fc_i+\frac{\partial c_i}{\partial t}\right)\omega_i+
\sum_{i=1}^sc_i \widetilde{\mathcal{L}}_f\omega_i,
\]

\noindent which belongs to span$\{\omega_1,\cdots,\omega_s,\widetilde{\mathcal{L}}_f\omega_1,\cdots,\widetilde{\mathcal{L}}_f\omega_s\}$
$\blacktriangleleft$}

\section{Computation of the controllable distribution}\label{AppendixDist}

We provide a property that plays a key role for the implementation of algorithm 4.
We denote by $\Delta_k$ the distribution returned by algorithm 4 at the $k^{th}$ step.
This property states that, at each step $k$, it suffices to compute the generators of $\Delta_k$ by performing simple operations (Lie derivative and time derivative) on the generators of $\Delta_{k-1}$.

\noindent We remind the reader the following property (that is the analogous that holds for algorithm 3):

\begin{pr}\label{PropertyConOld}
Let us consider a nonsingular distribution $\Delta$ spanned by the vector fields $f^1,\cdots,f^d$ and a smooth vector field $f$.
We have:

\[
\Delta\oplus[\Delta,~f]=\textnormal{span}\{f^1,\cdots,f^d,[f^1,~f],\cdots,[f^d,~f]\}
\]
\end{pr}

\proof{The reader is addressed to \cite{Isi95}, Remark 1.6.1
$\blacktriangleleft$}

\noindent We need to extend the above property to the case when, instead of the Lie bracket $[\cdot,~\cdot]$, we consider the new bracket defined in (\ref{EquationMixBracketTV}).

\noindent We have the following new property.

\begin{pr}\label{PropertyConNew}
Let us consider a nonsingular codistribution $\Delta$ spanned by the covectors $f^1,\cdots,f^d$ and a smooth vector field $f$. 
We have:

\[
\Delta\oplus\langle\Delta,~f\rangle=\textnormal{span}\{f^1,\cdots,f^d,\langle f^1,~f\rangle,\cdots,\langle f^d,~f\rangle\}
\]
\end{pr}

\proof{Obviously, span$\{f^1,\cdots,f^d,\langle f^1,~f\rangle ,\cdots,\langle f^d,~f\rangle \}$
is included in 
$\Delta\oplus\langle \Delta,~f\rangle $.

\noindent Let us prove the vice-versa. Let us consider a generic vector $v\in\Delta$. We have:

\[
v=\sum_{i=1}^dc_if^i
\]

\noindent We have:

\[
\langle v,~f\rangle =
\left[
\sum_{i=1}^dc_if^i, ~
f
\right]-
\frac{\partial}{\partial t} \sum_{i=1}^dc_if^i=
-\sum_{i=1}^d  \left(\mathcal{L}_fc_i+\frac{\partial c_i}{\partial t}\right)f^i+
\sum_{i=1}^dc_i \langle f^i,~f\rangle ,
\]

\noindent which belongs to span$\{f^1,\cdots,f^d,\langle f^1,~f\rangle ,\cdots,\langle f^d,~f\rangle \}$
$\blacktriangleleft$}

}

%\appendix
%\section{A summary of Latin grammar}    % Each appendix must have a short title.
%\section{Some Latin vocabulary}         % Sections and subsections are supported  
                                        % in the appendices.

\end{document}